\documentclass[11pt,dvips]{amsart}

\usepackage{amsmath}
\usepackage{amssymb}
\usepackage{amscd}
\usepackage{amsthm}
\usepackage{palatino}
\usepackage{euler}
\usepackage{latexsym}
\usepackage{epsfig}
\usepackage{graphics}
\usepackage{amsfonts}
\usepackage{psfrag}
\usepackage{graphicx}

\input xy
\xyoption{all}

\numberwithin{equation}{section}
\numberwithin{figure}{section}

\newtheorem{theorem}{Theorem}[section]
\newtheorem{lemma}[theorem]{Lemma}
\newtheorem{proposition}[theorem]{Proposition}

\newtheorem{definition}[theorem]{Definition}

\newtheorem{remark}[theorem]{Remark}

\newcommand{\C}{{\mathbb{C}}}
\renewcommand{\H}{{\mathbb{H}}}

\newcommand{\Z}{{\mathbb{Z}}}
\newcommand{\Q}{{\mathbb{Q}}}
\newcommand{\R}{{\mathbb{R}}}

\newcommand{\T}{{\mathbb{T}}}
\renewcommand{\P}{{\mathbb{P}}}

\newcommand{\algt}{\mathfrak{t}}

\newcommand{\iso}{\cong}
\newcommand{\into}{\hookrightarrow}
\newcommand{\onto}{\twoheadrightarrow}

\renewcommand{\mod}{/\!/}

\newcommand{\mmod}{{/\!\!/\!\!/\!\!/}}

\renewcommand{\ker}{{\operatorname k\!er}}
\newcommand{\Stab}{{\operatorname Stab}}
\newcommand{\Lie}{{\operatorname Lie}}
\newcommand{\image}{{\operatorname image}}
\newcommand{\Hom}{{\operatorname Hom}}
\newcommand{\age}{{\operatorname a\!ge}}

\newcommand{\id}{{\operatorname id}}
\newcommand{\im}{{\operatorname im}}

\newcommand{\hs}{{\hspace{3mm}}}
\newcommand{\hsm}{{\hspace{1mm}}}

\newcommand{\NHT}{NH^{*,\diamond}_T}

\newcommand{\NHTCot}{NH^{*,\diamond}_T(T^*\C^n)}
\newcommand{\NHGTCot}{NH^{*,\Gamma}_T(T^*\C^n)}

\begin{document}


\title{ORBIFOLD COHOMOLOGY OF HYPERTORIC VARIETIES}

\author{REBECCA F. GOLDIN}
\address{Mathematical Sciences MS 3F2\\ George Mason University\\
  4400 University Drive\\ Fairfax, VA 22030, USA}

\author{MEGUMI HARADA}
\address{Department of Mathematics and Statistics\\
McMaster University\\ 1280 Main Street West\\ Hamilton, ON L8S4K1, Canada}



\maketitle

\begin{abstract}
Hypertoric varieties are hyperk\"ahler analogues of toric varieties,
and are constructed as abelian hyperk\"ahler quotients $T^*\C^n \mmod
T$ of a quaternionic affine space.
Just as symplectic toric orbifolds are determined by labelled
polytopes, orbifold hypertoric varieties are intimately related to the
combinatorics of hyperplane arrangements.  By developing
hyperk\"ahler analogues of symplectic techniques developed by Goldin,
Holm, and Knutson, we give an explicit combinatorial description of
the Chen-Ruan orbifold cohomology of an orbifold hypertoric variety in
terms of the combinatorial data of a rational cooriented
weighted hyperplane arrangement ${\mathcal H}$.
We detail several
explicit examples, including some computations of orbifold Betti
numbers (and Euler characteristics).
\end{abstract}

\tableofcontents

\section{Introduction}\label{sec:intro}
The main result of this manuscript is an explicit combinatorial
computation of the Chen-Ruan orbifold cohomology $H^*_{CR}(M)$
\cite{CR04} of an orbifold hypertoric variety $M$. Hypertoric
varieties are hyperk\"ahler analogues of toric varieties, and were
first introduced by Bielawski and Dancer \cite{BieDan}, and
further studied by Konno \cite{Kon99}, \cite{Kon00} and
Hausel and Sturmfels \cite{HauStu}, among others.
Just as a symplectic toric orbifold is
determined by a labelled polytope, the theory of orbifold hypertoric
varieties is intimately related to the combinatorial data of a related
rational cooriented hyperplane arrangement ${\mathcal H}$. Our
description of $H^*_{CR}(M)$ is given purely in terms of this
arrangement ${\mathcal H}$. The fact that these hypertoric varieties
are constructed as hyperk\"ahler quotients $T^*\C^n \mmod T$ of a
quaternionic affine space \(\T^*\C^n \cong \H^n\) (via the
hyperk\"ahler analogue of the Delzant construction of K\"ahler toric
varieties) is crucial to our techniques.

Hyperk\"ahler quotients appear in many areas of mathematics. For
 instance, in representation theory, Nakajima's quiver varieties give
 rise to geometric models of representations (see e.g. \cite{Nak94},
 \cite{Nak98}, \cite{Nak01}). Furthermore, many moduli spaces appearing in physics,
 such as spaces of Yang-Mills instantons on $4$-manifolds or the
 solutions to the Yang-Mills-Higgs equations on a Riemann surface,
 arise via hyperk\"ahler quotient constructions. In each case, the
 study of topological invariants, such as cohomology rings or
 $K$-theory, of these quotients are of interest. In the case of
 hypertoric varieties, there are also close connections between the
 (ordinary or Borel-equivariant) cohomology rings of the varieties and
 the combinatorial theory of the corresponding hyperplane arrangements
 \cite{Kon99}, \cite{Kon00}, \cite{HauStu}, \cite{HP04}. Generalizing known such results
 to the orbifold case is of current
 interest. For example, recent work of Proudfoot and Webster
 \cite{ProWeb04}[Section 6] on the intersection cohomology of singular
 hypertoric varieties and the cohomology of their orbifold resolutions
 contains cohomological formulas which only apply in the unimodular
 case; it would be of interest to know whether there are orbifold
 versions of their statements.

In this paper, we focus on the combinatorics of the
hyperplane arrangement associated to
the {\bf Chen-Ruan orbifold cohomology} of orbifold hypertoric varieties.
Chen-Ruan orbifold cohomology rings were introduced in \cite{CR04} as
the degree $0$ piece of the Gromov-Witten theory of an orbifold,
following work in physics \cite{Zas96}. This ring carries,
in addition to the data of the usual singular cohomology ring of the
underlying space, more delicate information (e.g. about the orbifold
structure groups). Additively, $H^*_{CR}(M)$ is simply the usual
singular cohomology of the inertia orbifold $\widetilde{M}$ associated to
$M$; the product structure, on the other hand, is much more subtle,
incorporating the data of higher twisted sectors.
In this manuscript, we provide an
explicit presentation, via generators and relations, of this Chen-Ruan
cohomology ring for a class of orbifold hypertoric varieties.

Our approach is to develop hyperk\"ahler analogues
of the symplectic-geometric
techniques as introduced by Goldin, Holm, and Knutson in \cite{GHK05} to
compute Chen-Ruan orbifold cohomology. As in their work, we take advantage of the fact that a hypertoric
variety is by construction a global quotient of a manifold by a
torus.
We now briefly recall the main results of \cite{GHK05}.  Let $T$ be a
compact connected torus, let $N$ be a compact Hamiltonian $T$-manifold moment
map $\mu: N \to \algt^*$, and suppose that $\alpha$ is a regular value
of $\mu$. Then the inclusion \(\mu^{-1}(\alpha) \into N\)
induces a natural ring\footnote{In this paper, we take {\em rational}
  coefficients for all cohomology rings.}
homomorphism, often called the {\bf Kirwan map}:
\begin{equation}\label{eq:Kirwan}
\kappa: H^*_T(N) \onto H^*_T(\mu^{-1}(\alpha)) \cong H^*(N \mod_{\alpha}T),
\end{equation}
which is a surjection \cite{Kir84}. Here
$N \mod_{\alpha} T:= \mu^{-1}(\alpha)/T$ is by definition the
symplectic quotient of $N$ at $\alpha$.
The main result of \cite{GHK05} is an orbifold cohomology version of~\eqref{eq:Kirwan} for
abelian symplectic quotients. In other words, they show
that the inclusion \(\mu^{-1}(\alpha) \to N\) induces a surjective ring map (the ``{\bf orbifold
  Kirwan map}'')
\begin{equation}\label{eq:inertial-Kirwan}
\kappa_{NH}: NH^{*,\diamond}_T(N) \onto NH^{*,\diamond}_T(\mu^{-1}(\alpha)) \cong H^*_{CR}(M),
\end{equation}
where the domain is a new ring which they define: it is the
{\bf inertial cohomology ring} $NH_T^{*,\diamond}(N)$
of the $T$-space $N$.
(Here it is nontrivial that $\kappa_{NH}$ is a ring homomorphism; the same subtlety also arises
in the hyperk\"ahler case.) Moreover, they give an explicit description of the
kernel of $\kappa_{NH}$. Their proof relies on
symplectic-geometric properties of the fixed point sets $N^t$ for \(t
\in T\) in abelian Hamiltonian spaces, as well as on the original
Kirwan surjectivity result~\eqref{eq:Kirwan}.

In this paper we prove a parallel story in the hypertoric setting. A
direct hyperk\"ahler analogue of~\eqref{eq:inertial-Kirwan} is
nontrivial for several reasons, foremost among which is that it is not
known whether the hyperk\"ahler analogue of~\eqref{eq:Kirwan} is, in
general, surjective. However, in the special case where $M$ is a
smooth or orbifold hypertoric variety \cite{Kon00}, \cite{HauStu}, obtained
as a hyperk\"ahler quotient of $T^*\C^n$, a hyperk\"ahler analogue
of~\eqref{eq:Kirwan} does hold, allowing us to obtain results in this
setting.  The other nontrivial issues are the non-compactness of the
hypertoric varieties (in \cite{GHK05}, all orbifolds are assumed
compact) and the analysis of the hyperk\"ahler-geometric properties of
the fixed point sets $N^t \subseteq N$ for \(t \in T\). We deal with
these issues in Section~\ref{sec:surjection} to obtain the following.
Let $\mu_{HK}: M \to \algt^* \oplus \algt_{\C}^*$ denote the
hyperk\"ahler moment map on $T^*\C^n$, and $T^*\C^n \mmod T$ its
hyperk\"ahler quotient at a regular value $(\alpha, \alpha_\C) \in
\algt^* \oplus \algt_{\C}^*$ as described in
Section~\ref{sec:background-hypertoric}.

\begin{theorem}\label{th:surjectiontoCRcohomology}
Let $M$ be an orbifold
hypertoric variety $T^*\C^n \mmod T$.  There is a surjective ring homomorphism
\begin{equation}\label{eq:hk-inertial-Kirwan}
\kappa_{NH}^{\Gamma}: NH_T^{*,\Gamma}(T^*\C^n) \onto
H^*_{CR}(M),
\end{equation}
where $\Gamma$ is the subgroup of $T$ generated by finite stabilizers, $NH_T^{*,\Gamma}(T^*\C^n)$ is the $\Gamma$-subring of the
inertial cohomology ring $NH^{*,\diamond}_T(T^*\C^n)$, and $H^*_{CR}(M)$ is the Chen-Ruan cohomology of $M$.
\end{theorem}

The point of Theorem~\ref{th:surjectiontoCRcohomology} is that we can
in principle compute the orbifold cohomology of the hypertoric variety
$M = T^*\C^n \mmod T$ as a quotient of $NH^{*,\Gamma}_T(T^*\C^n)$ by
the kernel of~\eqref{eq:hk-inertial-Kirwan}. In the spirit of
\cite{Kon99,Kon00,HauStu,HP04}, we give an explicit algorithm for
computing both the domain $NH^{*,\Gamma}_T(T^*\C^n)$ and the ideal
\(\ker(\kappa_{NH}^{\Gamma})\) in terms of the combinatorics of a
central rational cooriented weighted hyperplane arrangement ${\mathcal
  H}_{cent}$ along with a choice of simple affinization ${\mathcal
  H}$.  This combinatorial data is obtained from the data of the
$T$-action on $T^*\C^n$ and an appropriate choice of level set of the
hyperk\"ahler moment map (explained in detail in
Section~\ref{sec:background-hypertoric}).
We now give a rough
statement of our main theorem, which gives a flavor of the ingredients
in the computation; the precise version is
Theorem~\ref{theorem:main-comb}.

\begin{theorem}\label{theorem:main-intro}
Let \(M = T^*\C^n \mmod T\) be an orbifold hypertoric variety.
Let ${\mathcal H} = \{H_i\}_{i=1}^n$ be a simple affine
rational cooriented hyperplane arrangement with positive normal
vectors $\{a_i\}_{i=1}^n$
associated to $M$ as described in
Section~\ref{sec:background-hypertoric}.
Then the Chen-Ruan cohomology of $M$
is given by
\[
H^*_{CR}(M)  \cong \Q[u_1, u_2, \ldots, u_n][\{\gamma_t\}_{t \in \Gamma}]
\bigg/ {\mathcal I} + {\mathcal J} + {\mathcal K}+\langle \gamma_{id}-1\rangle,
\]
where
\begin{itemize}
\item $\Gamma$ is a finite subgroup of $T$ determined by linear
independence relations among the $\{a_i\}_{i=1}^n$, made precise
in~\eqref{eq:intersection-kernels};
\item ${\mathcal I}$ is an ideal determined by $T$-weight data coming
from the action of $T$ on $T^*\C^n$ specified by ${\mathcal H}$, made
precise in Proposition~\ref{prop:inertial-mult};
\item ${\mathcal J}$ is an ideal generated by linear relations coming
from an exact sequence of Lie algebras \(0 \to \algt \to \algt^n \to
\algt^d \to 0\) given by the $T$-action on $T^*\C^n$ specified by
${\mathcal H}$, made precise in~\eqref{eq:Delzant-Lie};
\item ${\mathcal K}$ is determined by intersection data of the
  hyperplanes $H_i$ in ${\mathcal H}$, and made precise in
  Proposition~\ref{prop:kernel-kappa-t}.
\end{itemize}
\end{theorem}

In summary, this manuscript can be viewed in any of the following
ways. First, it is an example of an explicit computation of the
Chen-Ruan orbifold cohomology of hyperk\"ahler quotients, and a
further development, in the hyperk\"ahler setting, of the definition
and use of inertial cohomology as introduced in \cite{GHK05}. In
particular, we note that our methods would also apply to any class of
hyperk\"ahler quotients for which there exists an appropriate analogue
of the Kirwan surjection~\eqref{eq:Kirwan}. Similarly, although in
this manuscript we restrict our attention to $\Q$ coefficients for our
cohomology rings, if a $\Z$-coefficient analogue of the Kirwan
surjection for orbifold hypertoric varieties is proven, then our methods will
easily generalize to the setting of $\Z$ coefficients.  Second, it is
another exploration of the relationship between the geometry of
hypertoric varieties and the combinatorics of hyperplane arrangements.
Finally, it is the hyperk\"ahler-geometric analogue of the
algebraic-geometric description of the Chow ring of toric
Deligne-Mumford stacks in \cite{BCS05}.

In \cite{JT05}, Jiang and Tseng independently develop techniques for
an algebraic-geometric version of these results by defining
``hypertoric DM stacks'' using extended stacky fans, following work of
\cite{BCS05}. Their work applies to the sub-class of hypertoric
varieties M obtained by hyperk\"ahler quotients at regular values of
the form $(\alpha,0)$. In this case, there is a simple affine
hyperplane arrangement ${\mathcal H}$ in $\algt^*$ determined by the
data of a moment map for a residual torus action on $M$;  the results
of \cite{JT05} are phrased in terms of this arrangement ${\mathcal H}$.
Our results, on the other hand, apply to an orbifold hypertoric
variety obtained as a quotient at any regular value
$(\alpha,\alpha_{\C})$. This is because we do not
keep track of the hyperk\"ahler structure of the quotient (which does
depends on this choice of level set); the Chen-Ruan orbifold cohomology
of the quotient turns out to be {\em independent} of this choice,
i.e. is the same for any regular value. The
main difference between the approach taken in this manuscript and \cite{JT05}
is that Jiang and Tseng begin with the data of a
simple hyperplane arrangement ${\mathcal H}$ and then directly
construct the hypertoric DM stack associated to ${\mathcal H}$, which
has coarse moduli space the corresponding orbifold hypertoric variety.
As a result, they compute the product in the orbifold Chow ring
entirely in terms of the quotient hypertoric variety. In contrast, our
method is to work almost entirely ``upstairs'' on $T^*\C^n$ with a
linear $T$-action, before taking a hyperk\"ahler quotient. This  simplifies
some computations (as in \cite{GHK05}) by allowing us to work with
linear $T$-representations, and carries the information
of a family of hypertoric varieties at once.

Since orbifold Chen-Ruan cohomology reduces to ordinary cohomology
when $M$ is smooth, both our work and that of \cite{JT05} reduce to
the description of $H^*(M)$ given in \cite{Kon00}  (see also
\cite{HauStu}) in the case when $M$ is a smooth hypertoric variety. As
Jiang and Tseng illustrate \cite{JT05}, this can be useful to show
that the ordinary cohomology of a smooth hyperk\"ahler crepant
resolution of $\C^2/\Z_n$ as constructed by Kronheimer \cite{Kro89} is
isomorphic to the orbifold cohomology of $\C^2/\Z_n$, which can be
computed using \cite{FG03}.

We now give a summary of the contents of the paper.  In
Section~\ref{sec:background-hypertoric}, we give a brief account of
the construction of hypertoric varieties as a hyperk\"ahler quotient,
based on the data of a hyperplane arrangement.
In Section~\ref{sec:definition}, we briefly recall the
definition of inertial cohomology given in \cite{GHK05}. Then in
Section~\ref{sec:surjection}, we prove that there exists a surjection
in inertial cohomology as in~\eqref{eq:inertial-Kirwan}.
We give a
combinatorial description of the Chen-Ruan orbifold cohomology of a
hypertoric variety, based on the data of the hyperplane arrangement
${\mathcal H}$, in Section~\ref{sec:comb}.  In
Section~\ref{sec:example}, we work out in detail several explicit
examples, including some computations of orbifold Betti numbers and
orbifold Euler characteristics. The Appendix
(Section~\ref{sec:definitionChase}) contains a detailed discussion of
the isomorphism between inertial cohomology of a $T$-space $Z$ and
the Chen-Ruan cohomology of the quotient $X = Z/T$ (also discussed for
the compact case in \cite{GHK05}), as well as a
careful proof
of the correspondence between Chen and Ruan's definition of the
obstruction bundle with that used in the algebraic geometry literature
(e.g. \cite{FG03,BCS05}).

\section{Background: hypertoric varieties}\label{sec:background-hypertoric}

We first briefly describe the construction of hypertoric varieties in
order to set the notation and conventions to be used throughout the
rest of the paper. We refer the reader to \cite{BieDan,HP04,HauStu}
for a more leisurely account.

We begin with the hyperk\"ahler space $\H^n$,
 thought of as a holomorphic cotangent bundle \(T^*\C^n \cong
 \C^{2n}.\) This is a hyperk\"ahler manifold with real symplectic form
 $\omega_{\R}$ given by the identification with $T^*\C^n \cong
 \C^{2n}$ and $\omega_{\C}$ the canonical holomorphic symplectic
 form on a cotangent bundle. The standard linear diagonal action of the compact torus $T^n$
 on $\C^n$ induces an action on the holomorphic cotangent bundle
 $T^*\C^n$ which is hyperhamiltonian \cite{BieDan}. We will refer to
 this action as the {\bf standard hyperhamiltonian action of $T^n$ on
 $T^*\C^n$.} The hyperk\"ahler $T^n$-moment map $\tilde{\mu}_{HK} =
 (\tilde{\mu}_{\R}, \tilde{\mu}_{\C})$ on $T^*\C^n$ is given as
 follows. Let $\{u_i\}_{i=1}^n$ be a dual basis to $\{\varepsilon_i\}_{i=1}^n$ in $(\algt^n)^*$, and let
 \((z,w) = (z_1, \ldots, z_n, w_1, \ldots, w_n) \in T^*\C^n,\) where
 the $z_i$ are the base variables and the $w_i$ are the fiber
 variables. We have
\begin{align*}
\tilde{\mu}_{\R}(z,w) &= \frac{1}{2} \sum_{i=1}^n \left( \|z_i\|^2 - \|w_i\|^2 \right) u_i \in (\algt^n)^*,
\mbox{ and}\\
 \tilde{\mu}_{\C}(z,w) &= \sum_{i=1}^n z_i w_i u_i \in (\algt^n_{\C})^*.
\end{align*}

Let ${\mathcal H}_{cent} = \{H_i^{cent}\}_{i=1}^n$ be a central rational cooriented weighted hyperplane
arrangement  in
$(\algt^d)^*$ with positive normal vectors $\{a_i\}_{i=1}^n$ in
$\algt^d_{\Z}$. Here, ``weighted'' means that we do not require the
$a_i$ to be primitive vectors.  We now use this data to restrict the
$T^n$ action to that of a subtorus.  Let $\{\varepsilon_i\}_{i=1}^n$
be a basis of $\algt^n$. Define a linear map \(\beta: \algt^n \to
\algt^d\) by \(\beta(\varepsilon_i) = a_i \in \algt^d_{\Z}.\) Let
\(\algt = \algt^k := \ker(\beta) \subseteq \algt^n,\) where \(k =
n-d,\) with inclusion \(\iota: \algt^k \into \algt^n.\) This yields an
exact sequence
\begin{equation}\label{eq:Delzant-Lie}
\xymatrix @R=-0.2pc {
0 \ar[r] & \algt = \algt^k \ar[r]^{\iota} & \algt^n \ar[r]^{\beta} & \algt^d \ar[r] & 0, \\
         &                                & \varepsilon_i \ar@{|->}[r] & a_i & \\
}
\end{equation}
which on the one hand exponentiates to an exact sequence
\begin{equation}\label{eq:Delzant-Lie-group}
\xymatrix{
1 \ar[r] & T = T^k \ar[r]^{\exp\iota} & T^n \ar[r]^{\exp \beta} & T^d \ar[r] & 1,}
\end{equation}
and on the other hand
dualizes to the exact sequence
\begin{equation}\label{eq:Delzant-Lie-dual}
\xymatrix @R=-0.2pc {
0 \ar[r] & (\algt^d)^* \ar[r]^{\beta^*} \ar[r] & (\algt^n)^* \ar[r]^{\iota^*} & \algt^* = (\algt^k)^* \ar[r] & 0.\\
     & & u_i \ar@{|->}[r] & \lambda_i := \iota^* u_i & \\
}
\end{equation}
We will always assume that the set of integer vectors
\(\{a_i\}_{i=1}^n\) spans $\algt^d$ over $\Z$, so that the kernel $T =
T^k := \ker(\exp \beta)$ is connected; this assumption is also made in
\cite{HauStu}.

Now we restrict the $T^n$-action on $T^*\C^n$ to the subtorus $T$.
Let \(\lambda_i := \iota^* u_i \in (\algt^k)^*_{\Z}\) as in~\eqref{eq:Delzant-Lie-dual}. Then
$\lambda_i$ is the $T$-weight defining the action of the subtorus $T$
on the $i$-th coordinate of $\C^n$. Let \(\exp \lambda_i\) denote the
corresponding element in \(\Hom(T, S^1).\) Since the action of $T$ on
$T^*\C^n$ is given by the natural lift of that on $\C^n$, we have that
for \(t \in T, (z,w) \in T^*\C^n,\)
\begin{align}
t \cdot (z,w) & =
((\exp \lambda_1)(t)z_1, \ldots, (\exp \lambda_n)(t) z_n, (\exp \lambda_1)(t)^{-1}w_1, \ldots, (\exp \lambda_n)(t)^{-1} w_n).
\end{align}
The moment maps for the hyperhamiltonian $T$-action on $T^*\C^n$ are
given by composing $\tilde{\mu}_{HK}$ with the linear projection
\(\iota^*: (\algt^n)^* \to \algt^*.\) Thus we obtain the formulas
\begin{equation}\label{eq:mu-for-T}
\mu_{\R}(z,w) =
\frac{1}{2} \sum_{i=1}^n \left(\|z_i\|^2 - \|w_i\|^2 \right) \lambda_i \in \algt^*,
\quad \mbox{and} \quad
\mu_{\C}(z,w) = \sum_{i=1}^n z_i w_i \lambda_i \in \algt_{\C}^*.
\end{equation}
We will assume throughout that \(\lambda_i \neq 0, \forall i.\)

In order to specify the hyperk\"ahler quotient, we now pick a regular value
\((\alpha, \alpha_{\C}) \in \algt^* \oplus \algt^*_{\C} \cong
(\algt^*)^3\) at which to reduce. Any element \(\alpha \in
\algt^*\) specifies an affinization ${\mathcal H} = \{H_i\}_{i=1}^n$
of ${\mathcal H}_{cent}$ via the
equations
\[
H_i := \{ x \in (\algt^d)^*:   \langle x, a_i \rangle  = \langle - \tilde{\alpha}, \varepsilon_i \rangle \},
\]
where $\langle \cdot, \cdot \rangle$ denotes the natural pairing of a
vector space and its dual, and $\tilde{\alpha} \in (\algt^n)^*$ is a
lift of $\alpha$, i.e. \(\iota^* \tilde{\alpha} = \alpha.\) (A
different choice of lift just translates the whole hyperplane arrangement
by a constant.)

In particular, the choice of parameter \((\alpha,\alpha_{\C}) \in
(\algt^*)^3\) corresponds to {\em three} separate choices of
affinization of the central arrangement ${\mathcal H}_{cent}$. Hence a
hypertoric variety is determined by the combinatorial data of a
central weighted arrangement ${\mathcal H}_{cent}$ and, additionally,
$3$ choices of affinization of ${\mathcal H}_{cent}$. However, in the
case of our computation, this data can be simplified considerably.
This is because a preimage $\mu_{HK}^{-1}(\alpha,\alpha_{\C})$ of a
regular value of $\mu_{HK}$ is $T^n$-equivariantly diffeomorphic to
the preimage $\mu_{HK}^{-1}(\alpha', \alpha'_{\C})$ of any other
regular value; this can be seen by an argument essentially equivalent
to the proof of \cite{HP04}[Lemma 2.1]. Hence the Chen-Ruan cohomology
of \(T^*\C^n \mmod_{(\alpha,\alpha_{\C})} T\) can be seen to be
isomorphic to that of \(T^*\C^n \mmod_{(\alpha', \alpha'_{\C})} T,\)
as will be discussed further in
Remark~\ref{remark:CR-indep-of-alphaC}. In other words, the Chen-Ruan
cohomology of an orbifold hypertoric variety is determined by the
original central arrangement ${\mathcal H}_{cent}$, and is independent
of these choices of affine structures given by regular values. In
practice, however, it is useful to pick a convenient affinization with
which to work. Namely, if \(\alpha \in \algt^*\) is chosen such that
the corresponding affinization of ${\mathcal H}_{cent}$ is {\bf
  simple}\footnote{A hyperplane arrangement is simple if any subset of
  $\ell$ hyperplanes intersect in codimension $\ell$.}, then
$(\alpha,0) \in \algt^* \oplus \algt^*_{\C}$ is a regular value
\cite{BieDan}[Theorem 3.3]. Here, since the last two parameters are
both $0$, only the first factor gives rise to a nontrivial
affinization of ${\mathcal H}_{cent}$. We will denote by ${\mathcal H}
= \{H_i\}_{i=1}^n$ this {\em simple} affine rational cooriented
weighted hyperplane arrangement obtained from the data of ${\mathcal
  H}_{cent}$ and an appropriate \(\alpha \in \algt^*.\)

\section{Background: inertial cohomology}\label{sec:definition}

We begin with a brief account of inertial cohomology as developed
in \cite{GHK05}, which gives us a model for computing the orbifold
cohomology of the hyperk\"ahler quotients constructed in Section~\ref{sec:background-hypertoric}.
Readers already familiar with the definition of orbifold
cohomology in the sense of Chen and Ruan will find
Section~\ref{subsec:inertial-def-smile} straightforward,
since the product on the inertial cohomology of a $T$-space $Z$ is
defined precisely to mimic the Chen-Ruan product in the case that the
quotient \(X = Z/T\) is an orbifold, where $T$ acts with finite
stabilizers on $Z$. The contribution of \cite{GHK05} is to notice that
in other cases of $T$-spaces (such as Hamiltonian $T$-spaces), the
product on inertial cohomology can be described in terms of a
(different) product defined in terms of fixed point data. This
(different) product, which is easier to compute, is briefly recalled
in Section~\ref{subsec:robust}; it will play a key role in our
computation of the Chen-Ruan orbifold cohomology of hypertoric varieties.

\subsection{Inertial cohomology and the $\smile$ product}\label{subsec:inertial-def-smile}

Let $N$ be a stably complex
$T$-space. For any \(t \in T,\) let \(N^t\) denote the $t$-fixed points.
Since $T$ is abelian, each $N^t$ is also a $T$-space.

\begin{definition}\label{def:inertial-cohomology}
The {\bf inertial cohomology} of the space $N$ is, as a
$H^*_T(pt)$-module, given by
\begin{equation}\label{eq:inertial}
NH^{*,\diamond}_T(N) := \bigoplus_{t \in T} H^*_T(N^t),
\end{equation}
where the sum indicates the $\diamond$ grading, i.e. \(NH^{*,t}_T(N)
:= H^*_T(N^t)\).
\end{definition}

The $*$ grading on the left hand side is a {\em real}-valued grading defined in
\cite{Ruan02} which is obtained from the $*$ grading on the right hand side by a
shift depending (in this case) on the $T$-action; see
\cite{GHK05}[Section 3] for a detailed discussion.

\begin{remark}
In the case of the orbifold hypertoric varieties under
consideration in this paper, it will turn out that the $*$ grading is integral and
always even.
\end{remark}

Although the
stably complex structure does not enter into the definition of
the inertial cohomology as an additive group, it is an essential
ingredient in the definition of its product structure, which we now
discuss; first, however, we
warn the reader that the definition of the product on
$NH_T^{*,\diamond}(N)$, which we denote by \(a \smile b,\) is {\em
not} necessary for understanding the statement of our main
Theorem~\ref{th:surjectiontoCRcohomology}, but is necessary for the
proof. We include a brief definition only for completeness, and refer the reader to
\cite{GHK05}[Section 3] for details.

To describe the product, we make use of
the top Chern class of the ``obstruction bundle", which is a vector
bundle over connected components of certain submanifolds of $N$. More
specifically, let $t_1,t_2\in T$,  let $H = \langle t_1, t_2
\rangle$ be the subgroup they generate, and $N^H$ the submanifold of points fixed by $H$. For any connected component $Y$ of $N^H$, the normal bundle \(\nu(Y, N)\) of $Y$ in $N$ is naturally
equipped with an $H$-action.
We may decompose \(\nu(Y,N)\) into isotypic components with respect
to the $H$-action:
$$
\nu(Y, N) = \bigoplus_{\lambda\in\hat{H}} I_\lambda,
$$ where $\hat{H}$ denotes the character group of $H$.

\begin{definition}\label{def:logweight}
Let $\lambda \in \hat{H}$ and $t \in H$. For any connected component $Y$ of $N^H$, we define the {\bf logweight of $t$ with respect to $\lambda$}, denoted $a_{\lambda}(t)$,
to be the real number in $[0,1)$
such that $\lambda(t) = e^{2\pi i
a_\lambda(t)}$.
\end{definition}

Note that for any elements $t_1, t_2\in H$ and for any connected component of $N^H$, the sum
$a_\lambda(t_1)+a_\lambda(t_2)+a_\lambda((t_1t_2)^{-1})$
must be $0,1$, or $2$.

\begin{definition}\label{def:obstructionbundle}
The {\bf obstruction bundle} is a vector bundle over each component
$Y$ of $N^{H}$ specified by
$$ E|_Y := \bigoplus_{\lambda\in \hat{H} \atop
a_\lambda(t_1)+a_\lambda(t_2)+a_\lambda((t_1t_2)^{-1})=2} I_\lambda,
$$ where $\nu(Y,N)= \bigoplus I_\lambda$.
We write $E\rightarrow N^H$ to denote the union over all connected
components. Note that the dimension may vary over components. The {\bf virtual
fundamental class} $\varepsilon\in H_T^*(N^H)$ is given by
$$
\varepsilon := \sum_{Y \in \pi_0(N^H)} e(E|_Y),
$$ where 
$e(E|_Y)$ is the $T$-equivariant Euler class of $E|_Y$, considered as
an element of $H_T^*(Y)$.
\end{definition}
Now let $e_i: N^H \rightarrow
N^{t_i}$ for $i=1,2$, and $\overline{e}_3:N^H \rightarrow N^{t_1t_2}$
denote the natural inclusions. These induce pullbacks
$e_i^*:H_T^*(N^{t_i})\rightarrow H_T^*(N^H)$ for $i=1,2$ and the
pushforward $(\overline{e}_3)_*: H_T^*(N^H)\rightarrow
H^*_T(N^{t_1t_2})$.  Let $a\in NH_T^{*,t_1}(N)$ and $b\in
NH_T^{*,t_2}(N)$ be homogeneous classes in $\diamond$. Then we define
the product $a\smile b\in NH_T^{*,t_1t_2}(N)$, a homogeneous class in
$\diamond$, to be
$$
a\smile b:= (\overline{e}_3)_*(e_1^*(a)\cdot e_2^*(b)\cdot \varepsilon) \in H^*_T(N^{t_1t_2}) = NH^{*,t_1t_2}_T(N) \subset NH^{*,\diamond}_T(N),
$$
where the product $\cdot$ on the right hand side is the usual product on $H^*_T(N^H)$. Extending linearly, the product is defined on any two classes \(a,b,\).

\begin{remark}\label{remark:Gamma-subring}
It follows immediately from the definition of the product that, for any subgroup $\Gamma$ of $T$, there is a subring 
$NH_T^{*,\Gamma}(N)$ given by
$$
NH_T^{*,\Gamma}(N) = \bigoplus_{t\in\Gamma}H_T^*(N^t).
$$ We call this ring the {\em $\Gamma$-subring} of $NH_T^{*,\diamond}(N)$.
\end{remark}

\begin{remark}\label{remark:CR-indep-of-alphaC}
It is straightforward to show from the definition of $NH^{*,\diamond}_T$ that if $Z, Z'$
are stably complex $T$-spaces equipped with locally free $T$-actions
and there exists a $T$-equivariant diffeomorphism \(\phi: Z \to Z',\)
then \(NH^{*,\diamond}_T(Z) \cong NH^{*,\diamond}_T(Z')\) as graded
rings. Together with the proof given in the Appendix that the inertial cohomology of $Z$
is isomorphic to the Chen-Ruan cohomology of the quotient orbifold,
this justifies the claim in Section~\ref{sec:background-hypertoric}, i.e.
$H^*_{CR}(M)$ is indeed independent of the choice of regular value
$(\alpha, \alpha_{\C}) \in \algt^* \oplus \algt_{\C}^*$. In
particular, we may restrict without loss of generality to the case
$(\alpha,0)$.
\end{remark}

\subsection{The product on $NH_T^{*,\diamond}(N)$ when $N$ is robustly equivariantly injective}\label{subsec:robust}

In this section, we give a different description of the $\smile$
product which will be easier to use for our computations.
The $T$-space $N$ is {\bf robustly equivariantly injective} if the
natural inclusion $i:N^T \hookrightarrow N^t$ induces an injection in
equivariant cohomology $$i^*_t:H_T^*(N^t)\rightarrow H_T^*(N^T),$$ for
all $t\in T$. When $N$ satisfies this property, the product structure
on $NH_T^{*,\diamond}(N)$ can be described in terms of fixed point
data and the local structure of the $T$-action near fixed points;
see \cite{GHK05}.  Robust equivariant injectivity is a
strong condition: for instance, if $T$ acts locally freely on $N$, then \(N^T =
\emptyset\) and $N$ certainly cannot be robustly equivariantly
injective. On the other hand, Hamiltonian $T$-spaces are an important
source of examples of robustly equivariantly injective $T$-spaces.

For any component $F$ of the fixed point set $N^T$, $T$ acts on the
normal bundle to $F$. This representation splits into isotypic
component under the action:
\[
\nu(F, N) = \bigoplus_{\lambda\in\hat{T}} I_\lambda.
\]
Let $a\in NH_T^{*,t_1}(N)$ and $b\in NH_T^{*,t_2}(N)$ be homogeneous
classes in $\diamond$. Then $i^*_{t_1}(a)$ and $i^*_{t_2}(b)$ are
classes in $H^*_T((N^{t_1})^T)$ and $H^*_T((N^{t_2})^T)$,
respectively. Both of these rings are identified naturally with
$H_T^*(N^T)$. We define
\[
i^*_{t_1}(a)\star i^*_{t_2}(b)|_F:= i^*_{t_1}(a)|_{F}\cdot
i^*_{t_2}(b)|_{F}\cdot \prod_{I_\lambda \subset \nu(F, N)}
e(I_\lambda)^{a_\lambda(t_1)+a_\lambda(t_2)- a_\lambda(t_1t_2)},
\]
where $e(I_\lambda)\in H_T^*(F)$ is the equivariant Euler class of
$I_{\lambda}$, and all the products on the right hand side are
computed using the usual product in $H_T^*(F)$. Note that the exponent
is either 0 (if $a_\lambda(t_1)+a_\lambda(t_2)<1$) or 1 (otherwise).
By taking a sum over
the connected components and by extending linearly, this defines a new
product, which we call the {\bf $\star$ product,} on the image of
$i^*$ in \(\oplus_{t \in T} H^*(N^T).\) When $N$ is robustly
equivariantly injective, the map $i^*$ is injective, so the $\star$
product uniquely defines a product on $NH_T^{*,\diamond}(N)$.
By abuse of notation, we denote this product also as \(a
\star b,\) for \(a,b \in NH_T^{*,\diamond}(N).\)

The crucial fact, proven in \cite{GHK05}, is that {\em these two product
  structures agree,} i.e.
\[
a\smile b = a\star b,
\]
when $N$ is robustly equivariantly injective. Hence in the robustly
equivariantly injective case we may, for the purposes of computation,
work exclusively with the $\star$ product.  Note that
$N=T^*\C^n$ equipped with the $T$-action described in
Section~\ref{sec:background-hypertoric} is a Hamiltonian $T$-space, and in
particular it is robustly equivariantly injective. We use this in
Sections~\ref{sec:comb} and \ref{sec:example} to simplify the
combinatorics.

\section{Surjection in inertial cohomology}\label{sec:surjection}

Let $M$ be a hypertoric variety as constructed in
Section~\ref{sec:background-hypertoric}. Let \(Z :=
\mu_{HK}^{-1}(\alpha,\alpha_{\C}) \subseteq T^*\C^n\) be the level set of the
hyperk\"ahler moment map such that \(M = Z/T.\) In this section, we
show that the map
\begin{equation}\label{eq:inducedinertialmap}
\kappa_{NH}:\NHT(T^*\C^n)\rightarrow \NHT(Z)
\end{equation}
induced by the inclusion $i: Z\hookrightarrow T^*\C^n$
is a surjective ring homomorphism. In the appendix we prove that the latter ring is isomorphic to the orbifold cohomology of $M$ as (graded) rings, thus completing the proof of
Theorem~\ref{th:surjectiontoCRcohomology}.
An explicit description of both the domain and the kernel of
$\kappa_{NH}$, provided in Section~\ref{sec:comb}, will yield a
combinatorial description of $H_{CR}^*(M)$.
For the rest of the section, we will be largely following the outline of the proof of the
symplectic case in \cite{GHK05}. However, there are several new
considerations in the hyperk\"ahler case, which we will discuss as
they arise.

We must first justify why the inertial cohomology of the level $Z$
is defined. For this, it suffices
to observe that the normal bundle to $Z$ is
trivial since $Z$ is the preimage of a regular value of $\mu_{HK}$.
Therefore, the complex bundle $T(T^*\C^n)$ is a stabilization of $TZ$,
so $Z$ is also a stably complex $T$-space and
$NH^{*,\diamond}_T(Z)$ is well-defined.

Now consider the individual maps on Borel-equivariant cohomology
\begin{equation}\label{eq:Kirwan-map-classical}
\kappa_{NH}^t: H^*_T((T^*\C^n)^t) \to H^*_T(Z^t)
\end{equation} induced by the inclusions $Z^t \into (T^*\C^n)^t$.
Then we define the map on inertial cohomology to be the direct sum
of the $\kappa_t$, i.e.
\begin{equation}\label{eq:Kirwan-map-inertial}
\kappa_{NH} := \bigoplus_{t
  \in T} \kappa_{NH}^t: \NHTCot \to \NHT(Z).
\end{equation}

This is {\em a priori} only a map of $H_T^*(p)$-modules, {\em not necessarily a ring
homomorphism}.  Indeed, given an inclusion of the $T$-fixed point set
$N^T\hookrightarrow N$ of a Hamiltonian $T$-space, the induced map
\(NH^{*,\diamond}_T(N) \to NH^{*,\diamond}_T(N^T)\) on inertial
cohomology does not, in general, preserve the product structure since
the obstruction bundles are all trivial for $NH_T^{*,\diamond}(N^T)$.
However, if a $T$-equivariant inclusion $\iota: P
\into N$ behaves well with respect to the fixed point sets $N^t$ for
all $t \in T,$ then the obstruction bundles from Definition~\ref{def:obstructionbundle}
also behave well, and the induced map on $NH^{*,\diamond}_T$ is in
fact a ring homomorphism. We quote the following \cite{GHK05}[Proposition 5.1].

\begin{proposition}\label{prop:ring-map} {\bf (Goldin-Holm-Knutson)}
Let $N$ be a stably complex $T$-space. Let \(\iota: P \into N\)
be a $T$-invariant inclusion and suppose also that $P$ is transverse to
any $N^t$, \(t \in T.\) Then the map induced by inclusion \(\iota^*:
NH^{*,\diamond}_T(N) \to NH^{*,\diamond}_T(P)\) is a ring
homomorphism.
\end{proposition}

Thus, in order to check that the map $\kappa_{NH}$ is a ring
homomorphism, it suffices to check that the level set
$Z$ is transverse to any $(T^*\C^n)^t$.
We have the following general computation.

\begin{lemma}\label{lemma:transverse}
Let $T$ be a compact torus, and let $W$ be a hyperhamiltonian
$T$-space with moment map \(\mu_{HK} = (\mu_1, \mu_2,\mu_3): W \to
(\algt^*)^3.\) Assume $(\alpha_1, \alpha_2, \alpha_3)$ is a regular
value of $\mu_{HK}$, and let $Z$ denote the level set
\(\mu_{HK}^{-1}(\alpha_1, \alpha_2,\alpha_3) \subseteq W.\) Then $Z$
is transverse to $W^t$ for any \(t\in T.\)
\end{lemma}

\begin{proof}
The statement holds trivially \(Z \cap W^t = \emptyset,\). We assume that \(Z \cap
W^t \neq \emptyset\) and that $W^t$ is
connected; otherwise, we do the argument component by component.
Let \(y \in Z \cap W^t,\) and let \(\imath_t:
W^t \into W\) denote the inclusion. Since $W^t$ is a fixed point set
of a Hamiltonian $T$-action with respect to each symplectic form
$\omega_i$, $W^t$
is itself a
hyperhamiltonian $T$-submanifold of $W$, with moment map
\(\imath_t^*\mu_{HK}.\)

Since $(\alpha_1, \alpha_2, \alpha_3)$ is regular, the Lie algebra
\(\Lie(\Stab(y))\subseteq \algt\) of the stabilizer of $y$ is
$0$. In order to prove the transversality, it suffices to prove
that \(d(\imath_t^*\mu_{HK})_y\) is surjective. Since $W^t$ is
K\"ahler, and $\Lie(\Stab(y)) = \{0\}$,
\(d(\imath_t^*\mu_i)_y |_{J_i(T_y(T \cdot y))}\) is surjective onto
$\algt^*$ for each $i$, where $T\cdot y$ is the $T$-orbit through $y$ and $T_y(T \cdot y)$ its tangent space at $y$. Moreover, since $W^t$ is hyperk\"ahler, the three subspaces
$J_i(T_y(T\cdot y))$ for \(i=1,2,3\) are mutually orthogonal. In order
to show that \(d(\imath_t^*\mu_{HK})_y\) is surjective onto
$(\algt^*)^3$, it suffices to show that
\(d(\imath_t^*\mu_i)(J_j(T_y(T \cdot y))) = 0\) for \(i \neq j.\)
Without loss of generality we take \(i=2, j=1.\) For any \(X, Y \in \algt,\)
\begin{eqnarray*}
\langle d(\imath_t^* \mu_2)(J_1X_y^{\sharp}), Y \rangle &= &\omega_2(Y_y^{\sharp}, J_1 X_y^{\sharp})
\quad \mbox{ by definition of a moment map } \\
 &= & - g(Y_y^{\sharp}, J_2 J_1 X_y^{\sharp})
 \quad \mbox{ compatibility between } g, \omega_2 \\
 &= &  g(Y_y^{\sharp}, J_3 X_y^{\sharp})
 \quad \mbox{ quaternionic relation between the } J_i \\
 &= & - \omega_3(Y_y^{\sharp}, X_y^{\sharp})
 \quad \mbox{ compatibility between } g, \omega_3 \\
 &= &0 \qquad \mbox{since $T_y(T \cdot y)$ is isotropic with respect to
   $\omega_3$.}
\end{eqnarray*}
Thus $d(\imath_t^* \mu_{HK})_y$ maps the span of the three
subspaces \(J_i(T_y(T\cdot y)) \subseteq T_yW^t\) surjectively onto
$T_{(\alpha_1, \alpha_2, \alpha_3)} (\algt^*)^3 \cong (\algt^*)^3$.
This implies the
level set $Z = \mu_{HK}^{-1}(\alpha_1, \alpha_2, \alpha_3)$ is
transverse to $W^t$.
\end{proof}

Proposition~\ref{prop:ring-map} together with
Lemma~\ref{lemma:transverse} proves the following general fact.

\begin{proposition}\label{prop:ring}
Let $T$ be a compact torus, and let $W$ be a hyperhamiltonian
$T$-space with moment map \(\mu_{HK} = (\mu_1, \mu_2, \mu_3): W \to
(\algt^*)^3.\) Assume $(\alpha_1, \alpha_2, \alpha_3)$ is a regular
value of $\mu_{HK}$, and let $Z$ denote the level set
\(\mu_{HK}^{-1}(\alpha_1, \alpha_2, \alpha_3) \subseteq W.\) Then the map on
inertial cohomology induced by the inclusion \(Z \into W,\)
\[
NH^{*,\diamond}_T(W) \to NH^{*,\diamond}_T(Z),
\]
is a ring homomorphism.
\end{proposition}

In particular, in our case of hypertoric varieties, the map
$\kappa_{NH}$ defined in~\eqref{eq:Kirwan-map-inertial} is a ring
homomorphism. Now it remains to show that $\kappa_{NH} = \oplus_{t \in
  T} \kappa_{NH}^t$ is
surjective. To do this, we show that
\begin{equation}\label{eq:surjection-piece}
\kappa_{NH}^t: H^*_{T}((T^*\C^n)^t) \onto H^*_{T}(Z^t)
\end{equation}
is surjective for each \(t \in T\). We begin with an analysis of
these $t$-fixed point sets $(T^*\C^n)^t$.
A direct calculation shows that
\begin{align*}\label{eq:t-fixed-points}
(T^*\C^n)^t & := \{ \hsm (z,w) \in
  T^*\C^n \hsm \mid \hsm  z_i = w_i = 0 \hsm \mbox{if}\hsm
  (\exp\lambda_i)(t) \neq 1\hsm \} \\
  & \cong T^*\C^{S(t)} \\
\end{align*}
is a quaternionic affine subspace of $T^*\C^n$, where
$S(t):=\{ i\in\{1, 2, \ldots, n\}: (\exp \lambda_i)(t) = 1\}\subseteq \{1, 2, \ldots, n\}$ and
\begin{equation}\label{eq:def-CS}
\C^{S(t)} := \{(z_1, z_2, \ldots, z_n) \in \C^n: z_i = 0 \hsm \mbox{if} \hsm i \not \in S(t)\}.
\end{equation}
In addition, $T^*\C^{S(t)}$ is also a hyperhamiltonian $T$-space with
moment map given by \(\iota_t^*\mu_{HK},\) where \(\iota_t:
(T^*\C^n)^t \into T^*\C^n\) denotes the inclusion. This computation
allows us to conclude that~\eqref{eq:surjection-piece} is the {\em ordinary}
Kirwan map for the hypertoric subvariety $T^*\C^{S(t)}
\mmod_{(\alpha,\alpha_{\C})} T = Z^t/T$ of $M$. Since these maps are known to be surjective
\cite{Kon00,HauStu}, we have just proven the following.

\begin{theorem}\label{theorem:surjection-to-level-set}
Let $T^*\C^n$ be a
hyperhamiltonian $T$-space given by restriction of the standard
hyperhamiltonian $T^n$-action on $T^*\C^n$, where the inclusion \(T
\into T^n\) is determined
as in Section $2$.
Then the map on inertial cohomology induced by the inclusion
\(Z = \mu_{HK}^{-1}(\alpha,\alpha_{\C}) \into T^*\C^n,\)
\[
\kappa_{NH}: NH^{*,\diamond}_T(T^*\C^n) \to NH^{*,\diamond}_T(Z),
\]
is a surjective ring homomorphism.
\end{theorem}

Inertial cohomology is a direct sum
over infinitely many elements \(t \in T,\) so
Theorem~\ref{theorem:surjection-to-level-set} is not at all amenable to
computation. However, Theorem~\ref{theorem:surjection-to-level-set} can be
substantially simplified for computational purposes (in particular, it
can be made finite).  We first
establish some terminology.  Suppose a torus $T$ acts on a space
$Y$. Suppose \(y \in Y\) and the stabilizer group \(\Stab(y) \subseteq
T\) is finite. Then we call \(\Stab(y)\) a {\bf finite stabilizer
group.} Similarly, given a finite stabilizer group $\Stab(y)$, we call
any element \(t \in \Stab(y)\) a {\bf finite stabilizer (element).}
We let $\Gamma$ denote the subgroup in $T$ generated by finite
stabilizers. In the case of a linear $T$-action on $T^*\C^n$, this is
a finite subgroup of $T$, since the
$T$-action is determined by a finite set of weights.

By Remark~\ref{remark:Gamma-subring},
\begin{equation}\label{eq:def-Gamma-subring}
\NHGTCot := \bigoplus_{t \in \Gamma} NH^{*,t}_T(T^*\C^n) \subset NH^{*,\diamond}_T(T^*\C^n),
\end{equation}
is a subring of $\NHTCot$. We call this the {\bf
$\Gamma$-subring}. We now show that
$\kappa_{NH}$ is still
surjective when restricted to the $\Gamma$-subring.

\begin{theorem}\label{theorem:surjectivity-Gamma-subring}
Let $T^*\C^n$ be a
hyperhamiltonian $T$-space given by restriction of the standard
hyperhamiltonian $T^n$-action on $T^*\C^n$, where the inclusion \(T
\into T^n\) is determined
as in Section $2$.
Let $\Gamma$ be the subgroup in $T$ generated by finite stabilizers. Then the map
on the $\Gamma$-subrings of inertial cohomology induced by the
inclusion \(Z = \mu_{HK}^{-1}(\alpha,\alpha_{\C}) \into T^*\C^n,\)
\[
\kappa_{NH}^{\Gamma} := \kappa_{NH}|_{\NHGTCot}: \NHGTCot \to
NH^{*,\Gamma}_T(Z) \cong NH^{*,\diamond}_T(Z),
\]
is a surjective ring homomorphism.
\end{theorem}

\begin{proof}
Since the level set $Z$
is the preimage of a regular value of
$\mu_{HK}$, $T$ acts locally freely on
$Z$. In particular, \(Z^t = \emptyset\) if \(t \not \in \Gamma.\)
Hence for \(t \not \in \Gamma,\) the map \(\kappa_{NH}^t:
H^*_T((T^*\C^n)^t) \to H^*_T(Z^t)\)
is automatically $0$. Hence $\kappa_{NH}^t$ does not contribute to the
image of $\kappa_{NH}$, and \(\image(\kappa_{NH}) = \image(\kappa_{NH}
|_{\NHGTCot}).\) In particular, $\kappa_{NH}|_{\NHGTCot}$ is still
surjective.
\end{proof}

By Theorem~\ref{theorem:surjectivity-Gamma-subring}, we may
restrict our attention to the $\Gamma$-subring and the restricted ring
map \(\kappa_{NH}^{\Gamma}.\) The only remaining step to
complete the proof of Theorem~\ref{th:surjectiontoCRcohomology} is the
isomorphism of the inertial cohomology \(NH^{*,\diamond}_T(Z) \cong
NH^{*,\Gamma}_T(Z)\) of the level set with the Chen-Ruan cohomology
$H^*_{CR}(M)$ of the quotient \(M = Z/T.\) This would then imply that,
in order to compute $H^*_{CR}(M)$, it would suffice to compute the
domain and kernel of $\kappa_{NH}^{\Gamma}$.  
This explicit computation
is done in Section~\ref{sec:comb}. The isomorphism
$NH^{*,\diamond}_T(Z) \cong H^*_{CR}(M)$ mentioned above, for a general
stably complex $T$-space with locally free $T$-action, is
discussed in the compact case in \cite{GHK05}[Section 4]; we place a
detailed proof and its connection to the obstruction bundle in the algebraic geometry literature in an Appendix
(Section~\ref{sec:definitionChase}).

\section{The combinatorial description of $H^*_{CR}(M)$}\label{sec:comb}

We now come to the main result of this manuscript. Using the
inertial cohomology surjectivity result of
Section~\ref{sec:surjection} and the identification of inertial
cohomology with Chen-Ruan orbifold cohomology in
Section~\ref{sec:definitionChase}, we give in this section an explicit
description of the Chen-Ruan orbifold cohomology of hypertoric
varieties in terms of the combinatorial data of the hyperplane
arrangement ${\mathcal H}$. Here ${\mathcal H} = \{H_i\}_{i=1}^n$ is a
choice of simple affinization of the central arrangement ${\mathcal
  H}_{cent}$ determining the hyperhamiltonian $T$-action on $T^*\C^n$,
as detailed in Section~\ref{sec:background-hypertoric}.  We then work
out several concrete examples in Section~\ref{sec:example}.


We begin by stating our main theorem; for this, we must first set some
notation. Let $\mathcal{H}$ be a simple rational cooriented weighted
hyperplane arrangement. Suppose that \(S \subseteq \{1, 2, \ldots,
n\}\) such that \(\{a_j\}_{j \in S^c}\) is linearly independent in
$\algt^d$, where the $a_j$ are the integer normal vectors to the
hyperplanes $H_i$. Then
\begin{equation}\label{eq:Gamma-definition}
\Gamma_S := \bigcap_{i \in S} \ker(\exp \lambda_i) \subseteq T
\end{equation}
is a finite group. Let $\Gamma$ be the finite subgroup in $T$
generated by all such $\Gamma_S$.
For an element \(t \in \Gamma,\)
we define
\begin{equation}\label{eq:def-St}
S(t) := \{i: (\exp \lambda_i)(t) = 1\} \subseteq \{1, 2, \ldots, n\}.
\end{equation}
Let \(t_1, t_2 \in \Gamma.\) We define the
following subsets of \(\{1, 2, \ldots, n\}:\)
\begin{align}
A(t_1, t_2) & := \{ \hsm i \in S(t_1)^c \cap S(t_2)^c \hsm | \hsm
a_{\lambda_i}(t_1 t_2) = 0 \hsm \}, \nonumber \\
B(t_1, t_2) & := \left\{ \hsm j \in S(t_1)^c \cap S(t_2)^c \hsm \bigg| \hsm 
\begin{tabular}{l}
\(a_{\lambda_i}(t_1 t_2) \neq 0,\) \\ \(a_{\lambda_i}(t_1) + a_{\lambda_i}(t_2) - a_{\lambda_i}(t_1 t_2) = 0\) \\
\end{tabular}  \right\}, \label{eq:def-ABC}\\
C(t_1, t_2) & := \left\{\hsm k \in S(t_1)^c \cap S(t_2)^c \hsm \bigg| \hsm
\begin{tabular}{l} 
\(a_{\lambda_i}(t_1 t_2) \neq 0, \) \\  \(a_{\lambda_i}(t_1) +
a_{\lambda_i}(t_2) - a_{\lambda_i}(t_1 t_2) = 1\) \\ \end{tabular} 
\right\}, \nonumber
\end{align}
where $a_{\lambda}(t)$ is the logweight defined in
Definition~\ref{def:logweight}. Note that these sets partition the set
of indices corresponding to lines with nontrivial action by $t_1$ and
$t_2$. With this notation in place, we may state our main theorem.


\begin{theorem}\label{theorem:main-comb}
Let $T^*\C^n$ be a
hyperhamiltonian $T$-space given by restriction of the standard
hyperhamiltonian $T^n$-action on $T^*\C^n$, where the inclusion \(T
\into T^n\) is determined by the combinatorial data of ${\mathcal H}_{cent}$
as in Section $2$. Let ${\mathcal H}$ be a simple affinization of
${\mathcal H}_{cent}$.
Then for any regular value $(\alpha,\alpha_\C)$ the Chen-Ruan orbifold cohomology of the hyperk\"ahler quotient
$M := T^*\C^n \mmod_{(\alpha,\alpha_{\C})} T$
is given by
\[
H^*_{CR}(M) \cong \Q[u_1, \ldots, u_n][\{\gamma_t\}_{t \in \Gamma}] \Big/ {\mathcal I} + {\mathcal J} + {\mathcal K}+\langle \gamma_{id}-1\rangle,
\]
where the ideals ${\mathcal I}, {\mathcal J}, {\mathcal K}$ are
defined as follows. First,
\begin{equation*}
{\mathcal I} = \left< \gamma_{t_1} \gamma_{t_2} -
(-1)^{\sigma_{t_1t_2}}\left( \prod_{i \in A(t_1, t_2)} u_i^2 \right) \left( \prod_{j\in\atop{\small{B(t_1, t_2) \sqcup C(t_1,t_2)}}} u_j \right) \gamma_{t_1t_2} \hs \bigg{\vert} \hs t_1, t_2 \in \Gamma \right>
\end{equation*}
where $\sigma_{t_1t_2}={|A(t_1,t_2)|+|B(t_1,t_2)|}$. Second,
\begin{equation*}
{\mathcal J} = \left< \im(\beta^*) \right>.
\end{equation*}
Finally,
\begin{equation*}
{\mathcal K} = \sum_{t \in \Gamma} \left< \gamma_t \cdot \prod_{i \in L_t}
  u_i \hs \bigg{\vert} \hs \bigcap_{i \in L_t} H_i \cap \bigcap_{j \in
  S(t)^c} H_j = \emptyset \right>,
\end{equation*}
where $L_t$ denotes a (possibly empty) subset of $S(t)$.
\end{theorem}


The proof of Theorem~\ref{theorem:main-comb} involves three
steps. First, we must show that the subgroup generated by finite stabilizeres defined in
Section~\ref{sec:surjection} is indeed the group $\Gamma$ generated by
the $\Gamma_S$ in~\eqref{eq:Gamma-definition} above.
above.  Second, we prove that the ideals ${\mathcal I}, {\mathcal
  J}$ above are
exactly the relations which yield the inertial cohomology
$\NHGTCot$. Finally, we show that the ideal ${\mathcal K}$
exactly corresponds to the kernel of the inertial Kirwan map
$\ker(\kappa_{NH}^{\Gamma})$.


We begin with the first step, i.e. a description of the finite
stabilizer group $\Gamma$ associated to the given $T$-action on
$T^*\C^n$. As a bonus, we also give an (easy to compute) description
of the global orbifold structure groups that arise in the quotient
hypertoric variety. Let $M= T^*\C^n\mmod_{(\alpha,\alpha_{\C})} T$,
where the $T$-action on $T^*\C^n$ is determined by ${\mathcal
  H}_{cent}$. We have the following.

\begin{proposition}\label{prop:finite-stabilizers}
Let $T^*\C^n$ be a
hyperhamiltonian $T$-space given by restriction of the standard
hyperhamiltonian $T^n$-action on $T^*\C^n$, where the inclusion \(T
\into T^n\) is determined by the combinatorial data of ${\mathcal H}_{cent}$
as in Section~\ref{sec:background-hypertoric}. Let $\{a_i\}_{i=1}^n$
be the positive normal vectors defining the hyperplanes in ${\mathcal
  H}_{cent}$
and \(\lambda_i := \iota^*u_i\) as
in~\eqref{eq:Delzant-Lie-dual}.
\begin{enumerate}
\item A subgroup of $T$ is a finite stabilizer subgroup of a subvariety of $T^*\C^n$ if and only if it is of the form
\begin{equation}\label{eq:intersection-kernels}
\Gamma_S := \bigcap_{i \in S} \ker(\exp \lambda_i) \subseteq T,
\end{equation}
where \(S \subseteq \{1, 2, \ldots, n\}\) is such that \(\{a_j\}_{j
  \in S^c}\) is linearly independent in $\algt^d$. In particular, the
  subvariety \(M_S := T^*\C^S \mmod_{(\alpha,\alpha_{\C})} T\) has global
  orbifold structure group $\Gamma_S$.
\item The subgroup $\Gamma_S$ in~\eqref{eq:intersection-kernels} is isomorphic to
  \[
\left( \mbox{span}_{\Q}\{a_j\}_{j\in S^c} \cap \algt^d_{\Z} \right) \big/
\mbox{span}_{\Z}\{a_j\}_{j\in S^c}.
  \]
\item Any finite stabilizer \(t \in T\) occurs in a $\Gamma_S$ for
 $S$ such that \(\{a_j\}_{j \in S^c}\) forms a basis of $\algt^d$.
\end{enumerate}
\end{proposition}
\begin{proof}
We begin with a general computation.
Let \((z,w) \in T^*\C^n.\) Recall that the action of the subtorus \(T \subseteq T^n\) is given by composing the homomorphism
\begin{equation}
T \to T^n, \quad t \mapsto ((\exp \lambda_1)(t), \ldots, (\exp \lambda_n)(t)),
\end{equation}
with the standard linear action of $T^n$ on $T^*\C^n$. It is immediate that
\begin{equation}\label{eq:stabilizer-v1}
\Stab(z,w) = \left\{ t \in T: \mbox{if either} \hsm z_i \neq 0 \hsm \mbox{or} \hsm w_i \neq 0, \hsm \mbox{then} \hsm (\exp \lambda_i)(t) = 1 \right\}.
\end{equation}
Now define \(S(z,w) := \left\{ i \hsm | \hsm
z_i \neq 0 \hsm \mbox{or} \hsm w_i \neq 0 \right\}.\)
Then~\eqref{eq:stabilizer-v1} becomes
\begin{equation}
\Stab(z,w) = \bigcap_{i \in S(z,w)} \ker(\exp \lambda_i).
\end{equation}
In particular, $\Stab(z,w)$ is finite if and only if the set
\(\{\lambda_i\}_{i \in S(z,w)}\) spans $\algt^*$, or equivalently, the intersection
\(\cap_{i \in S(z,w)} \ker(\lambda_i) = \{0\}.\)
By the exactness of the sequence~\eqref{eq:Delzant-Lie}, this is equivalent
to the condition that \(\{a_j\}_{j \in S(z,w)^c}\) is linearly
  independent in $\algt^d$.
Conversely, given a subset $S$ with $\{a_j\}_{j \in S^c}$ linear
independent, any $(z,w) \in T^*\C^n$
such that \(z_i = w_i = 0\) for \(i \not \in
S,\) and \(z_i \neq 0\) or \(w_i \neq 0\) for \(i \in S,\) will
have stabilizer exactly \(\cap_{i \in S} \ker(\exp \lambda_i).\)
Moreover, the argument above immediately
implies that $M_S$ has global orbifold structure group
$\Gamma_S$. This proves the first claim.

To prove the second claim, we will produce a map $\varphi$ from
\(\left( \mbox{span}_{\Q}\{a_j\}_{j\in S^c} \cap \algt^d_{\Z}
\right)\) to $T$, which we will show takes values in $\Gamma_S$. For the remainder of
this computation, we identify $T$ with \(\ker(\beta)/(\ker(\beta) \cap
\algt^n_{\Z})\) by~\eqref{eq:Delzant-Lie} and $S^1$ with $\R/\Z$.  In
this language, \([X] \in T\) is in $\Gamma_S$
exactly when any representative \(X = \sum_ix_i\varepsilon_i \in
\ker(\beta)\) of $[X]$ has the property that \(x_i \in \Z\) for all
\(i \in S.\)
We begin by constructing the map $\varphi$. Let \(y \in
\mbox{span}_{\Q}\{a_j\}_{j\in S^c} \cap \algt^d_{\Z}.\) Since
\(\mbox{span}_{\Z}\{a_i\}_{i=1}^n =
\algt^d_{\Z}\) by assumption, there exist linear combinations
\begin{equation}\label{eq:Z-and-Q-combination}
y = \sum_{i=1}^n c_i a_i, \quad
\mbox{and} \quad y = \sum_{j \in S^c} d_j a_j,
\end{equation}
where \(c_i \in \Z, d_j \in \Q,\) and the second linear
combination is unique.
Let \(X = \sum_i x_i \varepsilon_i \in \algt^n\) where
\[
x_k := \left\{ \begin{array}{ll} c_k & \mbox{if $k \in S$} \\
                                c_k - d_k & \mbox{if $k \in S^c$.} \\
              \end{array}
      \right.
\]
Then by construction $x$ represents an element
in $\Gamma_S$, and we define \(\varphi(y) := [x]
\in T.\) A different choice of $\Z$-linear combination
in~\eqref{eq:Z-and-Q-combination} yields the same $[x]$, so
$\varphi$ is well-defined.
Furthermore, by definition, if \(y \in \mbox{span}_{\Z}\{a_j\}_{j \in
  S^c}, \varphi(y)\) is trivial in $T$, so
$\varphi$ also factors through the quotient
\[
\left( \mbox{span}_{\Q}\{a_j\}_{j\in S^c} \cap \algt^d_{\Z} \right) \big/
\mbox{span}_{\Z}\{a_j\}_{j\in S^c}.
\]
The map $\varphi$ preserves additive structures, hence is a
 homomorphism. Furthermore, $\varphi$ is
an injection since if \(\varphi(y) \in \ker(\beta) \cap
\algt^n_{\Z},\) then the coefficients \(d_j\)
in~\eqref{eq:Z-and-Q-combination} are integers, and hence \(y \in
\mbox{span}_{\Z}\{a_j\}_{j \in S^c}.\) Finally, to see that $\varphi$
is surjective, let \(X \in \ker(\beta)\) be a representative for an
element in $\Gamma_S$, with coordinates \(c_i\)
for \(i \in S,\) \(x_j\) for \(j \in S^c.\)
Then
\(y := \sum_{i \in S} c_i a_i\) has the property that \(\varphi(y) = [X],\) so $\varphi$ is surjective.
  Hence $\varphi$ is an isomorphism, as desired.

Finally, since we always have
\[
\bigcap_{i \in S'} \ker(\exp \lambda_i) \subseteq \bigcap_{i \in S} \ker(\exp \lambda_i)
\]
for any \(S \subseteq S',\) in order to identify the finite stabilizer
elements in $T$, it suffices to consider the {\em minimal} subsets $S$
such that \(\cap_{i\in S} \ker \lambda_i = \{0\},\) or equivalently,
{\em maximal} linearly independent sets \(\{a_j\}_{j \in S^c},\)
i.e. bases of $\algt^d$. This proves the final claim.
\end{proof}

Thus, in order to compute $\Gamma$, it suffices to find the subsets
$\{a_j\}_{j \in S^c}$ in \(\{a_i\}_{i=1}^n\) which form a basis of
$\algt^d$. We also note that the subvarieties $M_S$ map under the moment map for $M$
 to the intersection of the hyperplanes $\bigcap_{j \in S^c}
H_j$, so can easily be identified in the combinatorial picture using
${\mathcal H}$.

We now proceed to the second step, i.e. we describe the product structure on $\NHGTCot$.

\begin{proposition}\label{prop:inertial-mult}
Let $T^*\C^n$ be a hyperhamiltonian $T$-space given by restriction of
the standard hyperhamiltonian $T^n$-action on $T^*\C^n$.  Let $\NHGTCot$ be the
$\Gamma$-subring of the inertial cohomology ring $\NHTCot$, and let
\(\lambda_i := \iota^* u_i\)
as in~\eqref{eq:Delzant-Lie-dual}.  Then, as a
graded $H^*_T(pt; \Q)$-algebra,
\[
\NHGTCot \cong \Q[u_1,u_2, \ldots, u_n][\{\gamma_t\}_{t \in \Gamma}] \bigg/
     {\mathcal I}+{\mathcal J} + \langle \gamma_{\id} -1\rangle,
\]
where the ideal ${\mathcal I}$ is generated by the relations
\begin{equation}\label{eq:inertial-relations}
\gamma_{t_1} \gamma_{t_2} = (-1)^{|A(t_1,t_2)|+|B(t_1,t_2)|}\left( \prod_{i \in A(t_1, t_2)} u_i^2 \right) \left( \prod_{j \in B(t_1, t_2) \sqcup C(t_1,t_2)} u_j \right)\gamma_{t_1 t_2} ,
\end{equation}
with the sets $A(t_1,t_2),B(t_1,t_2),C(t_1,t_2)$ as defined
in~\eqref{eq:def-ABC}, and ${\mathcal J}=\langle\im(\beta^*)\rangle$.
\end{proposition}

\begin{remark}The grading is given by $\deg u_i=2$ for all $i$, and $\deg \gamma_t = 2age(t)$, as specified in the Appendix.
\end{remark}
\begin{proof}
Recall that the $\Gamma$-subring $\NHGTCot$ is by definition given by
\[
\NHGTCot := \bigoplus_{t \in \Gamma} H^*_T((T^*\C^n)^t).
\]
Since each $t$-th graded piece is the $T$-equivariant cohomology of a
contractible space, it has a single generator as a $H^*_T(pt)$-module.
Let $\gamma_t$ denote the element in $\NHGTCot$ which is equal to $0$ for
each $h$-graded piece with \(h \neq t, h \in \Gamma,\) and which is
equal to the generator \(1 \in H^*_T((T^*\C^n)^t) \cong H^*_T(pt)\) in
the $t$-th graded piece. Then $\NHGTCot$ is generated as a
$H^*_T(pt)$-module by these $\{\gamma_t\}_{t\in \Gamma}$.  Hence
in order to determine the multiplicative structure, it suffices to
find the product relations among these generators $\gamma_t, t \in
\Gamma$. Also, by the exact sequence~\eqref{eq:Delzant-Lie}, we may
identify
\[
H^*_T(pt;\Q) \cong H^*_{T^n}(pt;\Q) \big/ {\mathcal J} \cong \Q[u_1,
  \ldots, u_n] \big/ {\mathcal J}.
\]

Since $T^*\C^n$ is robustly equivariantly injective, we may compute
all products in terms of the $\star$ product instead of the $\smile$
product, as was explained in Section~\ref{subsec:robust}. By our
assumptions on ${\mathcal H}$, all the weights $\lambda_i$ defining
the action of $T$ on $T^*\C^n$ are non-zero, and hence the only
$T$-fixed point is the origin $\{0\} \in T^*\C^n$. The $T$-weights of
the action on the normal bundle $\nu(F, N)$ to $F := (T^*\C^n)^T =
\{0\}$ are the $2n$ weights \(\{\pm \lambda_i\}_{i=1}^n.\) Since the
weights come in pairs, the definition of the $\star$ product yields
\begin{equation*}\label{eq:inertial-product-relation}
\gamma_{t_1} \star \gamma_{t_2} = \gamma_{t_1 t_2} \cdot \prod_{i=1}^n \left( (\lambda_i)^{a_{\lambda_i}(t_1) + a_{\lambda_i}(t_2) - a_{\lambda_i}(t_1t_2)} \cdot  (-\lambda_i)^{a_{-\lambda_i}(t_1) + a_{-\lambda_i}(t_2) - a_{-\lambda_i}(t_1t_2)} \right).
\end{equation*}
If $i\in S(t_1)\cup S(t_2)$, then either \(a_{\lambda_i}(t_1) = 0\) or
 \(a_{\lambda_i}(t_2) = 0,\) and the corresponding exponent is
 0. Suppose \(i \in S(t_1)^c \cap S(t_2)^c\). We now take
 cases. Suppose that \(i \in A(t_1, t_2).\) In this case,
 \(a_{-\lambda_i}(t_{\ell}) = 1 - a_{\lambda_i}(t_{\ell})\) for either
 \(\ell = 1\) or $2$, $a_{\lambda_i}(t_1)+a_{\lambda_i}(t_2)=1$ and the $i$-th term in the
 product above is
 $-\lambda_i^2$. Similar computations show that if \(j \in B(t_1,
 t_2),\) then the $j$-th term is equal to $-\lambda_j$, and if \(k \in
 C(t_1, t_2),\) then the $k$-th term is $\lambda_k$. Finally, given
 the identification of $H^*_T(pt;\Q)$ with $H^*_{T^n}(pt;\Q)/{\mathcal
 J}$, a representative of $\lambda_i \in H^2_T(pt;\Q)$ is given by
 $u_i \in H^2_{T^n}(pt;\Q)$. 
\end{proof}

The third and final step is to determine the kernel of the inertial Kirwan map $\kappa^{\Gamma}_{NH}$. Since $\kappa^{\Gamma}_{NH} = \oplus_{t \in \Gamma} \kappa_{NH}^t$ is a direct sum of maps \(\kappa^t_{NH}: NH^{*,t}_T((T^*\C^n)^t) \to NH^{*,t}_T(Z^t)\) for \(t \in \Gamma,\) it suffices to compute the kernel of each $\kappa^t_{NH}$ separately.

Suppose \(t \in \Gamma.\) Then, as observed in
Section~\ref{sec:surjection}, $\kappa_{NH}^t$ is the Kirwan map in
usual cohomology $H^*(-;\Q)$ for the hyperk\"ahler Delzant
construction of an orbifold hypertoric variety. This map is known to
be surjective and the kernel has been explicitly computed. We quote
the following.


\begin{theorem}\label{theorem:Hausel-Sturmfels} {\bf (Hausel-Sturmfels)}
Let $T^*\C^n$ be a hyperhamiltonian
$T$-space given by restriction of the standard linear hyperhamiltonian
$T^n$-action on $T^*\C^n$, where the inclusion \(T \into T^n\) is
determined by the combinatorial data of ${\mathcal H}_{cent}$ as in
Section~\ref{sec:background-hypertoric} and let ${\mathcal H}$ be a
simple affinization of ${\mathcal H}_{cent}$.
Then the ordinary cohomology ring of the
orbifold hypertoric variety $M = T^*\C^n \mmod_{(\alpha,\alpha_{\C})} T$ is
given by
\begin{equation}\label{eq:Hausel-Sturmfels}
H^*(M;\Q) \iso \Q[u_1, \ldots, u_n] \bigg/ {\mathcal J} + \left< \prod_{i\in L} u_i
\hsm \bigg{\vert} \hsm \bigcap_{i \in L} H_i = \emptyset \right>,
\end{equation}
where $L$ denotes a subset of $\{1, 2, \ldots, n\}$.
\end{theorem}


We will apply Theorem~\ref{theorem:Hausel-Sturmfels} to the hyperhamiltonian action of $T$ on $T^*\C^{S(t)}$ for each
\(t \in \Gamma.\) Since the expression in~\eqref{eq:Hausel-Sturmfels} uses the combinatorial data of
the hyperplane arrangement corresponding to this action, our first task will be to describe explicitly the arrangement
${\mathcal H}_t$ for each \(t \in \Gamma\) in terms of the original ${\mathcal H}$.
Since the $T$-action is defined as a
restriction of the original $T$-action on $T^*\C^n$, the action on
$T^*\C^{S(t)}$ is given by the composition \(\iota_t := \pi_t \circ \iota,\)
\[
\xymatrix{
\algt \ar[r]^{\iota} \ar @/^2pc/[rr]^{\iota_t} & \algt^n \ar[r]^{\pi_t} & \algt^{S(t)} \\
}
\]
where $\iota$ is the inclusion coming from the original Delzant
sequence~\eqref{eq:Delzant-Lie}, and $\pi_t$ is the natural
projection to the subspace, given by \(\varepsilon_i \mapsto [\varepsilon_i], \forall i \in
S(t), \varepsilon_j \mapsto 0, \forall j \in S(t)^c.\)
A simple linear algebra argument
together with the commuting diagram
\begin{equation}\label{eq:comm-diag}
\xymatrix{
0 \ar[r] & \algt \ar[r]^{\iota}
\ar[dr]^{\iota_t} & \algt^n \ar[r]^{\beta} \ar[d]^{\pi_t} & \algt^{d=n-k} \ar[d] \\
 & & \algt^{S(t)} \cong \algt^n \big/ {\operatorname span}\left< \varepsilon_i
\right>_{i \in S(t)^c} \ar[r]^(0.55){\beta_t} & \algt^d \big/
{\operatorname span}\left< a_i \right>_{i \in S(t)^c} \\
}
\end{equation}
(where the top exact sequence is that in~\eqref{eq:Delzant-Lie}, the
right vertical arrow is the natural projection, and $\beta_t$ is the
composition of $\beta$ with the natural projection) shows that the map
$\iota_t$ fits into the exact sequence
\begin{equation}\label{eq:Delzant-Lie-t}
\xymatrix{
0 \ar[r] & \algt \ar[r]^{\iota_t} & \algt^{S(t)} \ar[r]^(.3){\beta_t} &
\algt^d \big/{\operatorname span}\left< a_i \right>_{i \in S(t)^c} \ar[r] & 0.
}
\end{equation}
From this sequence we will be able to deduce the structure of the arrangement ${\mathcal H}_t$
We note, however, the
sequence~\eqref{eq:Delzant-Lie-t} is not necessarily a standard
Delzant exact sequence as in
Section~\ref{sec:background-hypertoric}. This is because it is
possible to have \(\beta_t([\varepsilon_i]) = 0\) for some \(i \in
S(t),\) whereas this does not occur for a standard Delzant
construction. This poses no serious problems, as will be discussed in
more detail later. The relevant
combinatorial data must therefore be contained in the {\em non-zero}
images, $[a_i] \neq 0, i \in S(t)$.

Using~\eqref{eq:Delzant-Lie-t}, we may now explicitly describe the
hyperplane arrangement ${\mathcal H}_t$ for $t \in \Gamma$.  First,
${\mathcal H}_t$ sits naturally in the dual of the Lie algebra
\(\algt^d \big/{\operatorname span}\left< a_j \right>_{j \in
  S(t)^c},\) which is a subspace of
$(\algt^d)^*$. Specifically, it is the annihilator of the
$|S(t)^c|$-dimensional subspace \({\operatorname
  span}\left<a_j\right>_{j \in S(t)^c}.\) Hence, up to an affine
translation, it may be identified with the intersection
\begin{equation}\label{eq:Ht-space}
\bigcap_{j \in S(t)^c} H_j \subseteq (\algt^d)^*,
\end{equation}
where $H_j$ is the hyperplane orthogonal to $a_j$ in the original hyperplane arrangement ${\mathcal H}$.
Moreover, by analyzing the dimensions of $T^d$-orbits in the subvariety $M_{S(t)} \subseteq M$, it is straightforward to see that the affine hyperplanes in ${\mathcal H}_t$ are exactly given by the intersections
\begin{equation}\label{eq:Ht-hyperplanes}
H_i \cap \bigcap_{j \in S(t)^c} H_j  \subseteq \bigcap_{j \in S(t)^c} H_j,
\end{equation}
where \(i \in S(t)\) such that \(a_i \not \in \mathrm{span} \left< a_j \right>_{j \in S(t)^c}.\)

The only indices \(i
\in S(t)\) which give non-empty hyperplanes in ${\mathcal H}_t$ are
those for which \([a_i] = \beta([\varepsilon_i]) \neq 0.\) Thus, in
addition to the standard Delzant construction for ${\mathcal H}_t$
with certain basis vectors $[\varepsilon_i]$ mapping to the
corresponding normal vector $[a_i] \neq 0$ defining a hyperplane in
${\mathcal H}_t$, we have in this case also some extra basis
elements corresponding to the $i \in S(t)$ for which
\(\beta[\varepsilon_i]=0.\) This poses no problems, because any such
extra indices correspond to a subtorus of $T$ acting {\em standardly}
on a quaternionic affine space, the moment map for which has level
sets precisely equal to group orbits and hence has trivial
hyperk\"ahler quotient. In particular, the addition of such extra
indices leaves the corresponding hypertoric variety topologically
unchanged, allowing us (with only slight modifications) to use the
known theorems for hypertoric varieties built via a standard Delzant
construction.

We have the following Proposition.

\begin{proposition}\label{prop:kernel-kappa-t}
Let $T^*\C^n$ be a
hyperhamiltonian $T$-space given by restriction of the standard
hyperhamiltonian $T^n$-action on $T^*\C^n$, where the inclusion \(T
\into T^n\) is determined by the combinatorial data of ${\mathcal H}_{cent}$
as in Section~\ref{sec:background-hypertoric}, and let ${\mathcal H}$
be a simple affinization of ${\mathcal H}_{cent}$.
Let \(t \in \Gamma\). Then the ordinary
cohomology of the hypertoric subvariety \(M_S = T^*\C^{S(t)} \mmod T\)
is given by
\[
H^*(M_S;\Q) \cong \Q[u_1, u_2, \ldots, u_n] \bigg/ {\mathcal J} + {\mathcal K}_t,
\]
where
\begin{equation}\label{eq:def-ideal-Kt}
{\mathcal K}_t = \left< \prod_{i \in L_t} u_i \hsm \bigg{\vert} \hsm \bigcap_{i \in L_t} H_i \cap \bigcap_{j \in S(t)^c}H_j = \emptyset \right>,
\end{equation}
where $L_t$ denotes a (possibly empty) subset of $S(t)$, and ${\mathcal J} = \left< \im(\beta)^* \right>.$
\end{proposition}

\begin{proof}
We begin by observing that the domain $H^*_T(T^*\C^{S(t)}) \cong H^*_T(pt)$ of the map $\kappa_{NH}^t$ may also be identified as \(H^*_{T^n}(pt) \big/ \left<\im(\beta^*)\right> \cong \Q[u_1, u_2, \ldots, u_n] \big/\left< \im(\beta^*) \right>,\) by~\eqref{eq:Delzant-Lie}.

A subtlety that arises here is the presence of global stabilizers for
the $T$-action on the subsets $(T^*\C^n)^t = T^*\C^{S(t)}$ for $t \in \Gamma$ a
non-trivial finite stabilizer. Clearly, if $t \neq \id$, then by
definition $(T^*\C^n)^t$ has some non-trivial global stabilizer
$\Gamma_t$. Hence the $T$-action on $T^*\C^{S(t)}$ is not effective
and in particular does not arise from a standard Delzant construction
(since any such is effective).
However, since $\Gamma_{t}$ is finite, $T/\Gamma_{t}$ is again a torus of
dimension $\dim(T)$, and the inclusion maps on the level of Lie
algebras are identical. The same holds at the level of cohomology
rings, and hence the computation with global finite stabilizer is
identical to the computation in the usual hyperk\"ahler Delzant
construction. (Put another way, the essential data for the computation
is in the maps on Lie {\em algebras}.)

Putting together the description given in~\eqref{eq:Ht-space} and
\eqref{eq:Ht-hyperplanes} of the hyperplane arrangement ${\mathcal
H}_t$ associated to~\eqref{eq:Delzant-Lie-t},
Theorem~\ref{theorem:Hausel-Sturmfels}, and the commutative
diagram~\eqref{eq:comm-diag}, we see that the kernel of
$\kappa_{NH}^t$ is generated by the relations given
in~\eqref{eq:def-ideal-Kt}. Note that if \(i \in S(t)\) such that
\(\beta[\varepsilon_i] = 0,\) then
\[
H_i \cap \bigcap_{j \in S(t)^c} H_j = \emptyset,
\]
so \(u_i \in {\mathcal K}_t\).
Finally, we observe that if \(\cap_{j \in S(t)^c} H_j = \emptyset\)
(for instance if \(t\in \Gamma\) is not a finite stabilizer) then
\(\kappa^t_{NH} \equiv 0\) since in this case \(Z^t = \emptyset.\)
This in particular implies that we must have \(\gamma_t \in {\mathcal
K}_t,\) which is implied by our convention that we can take \(L_t =
\emptyset\) in the relations above and hence \(\gamma_t \cdot 1 = \gamma_t
\in {\mathcal K}_t.\) The result follows.
\end{proof}

We may now prove our main theorem, which gives a full combinatorial description of the Chen-Ruan orbifold cohomology of the orbifold hypertoric variety $M$.

\begin{proof} {\bf (Proof of Theorem~\ref{theorem:main-comb})}
We will give a description of the Chen-Ruan orbifold cohomology of $M$
as a quotient of $NH^{*,\Gamma}_T(T^*\C^n)$ by the kernel of
$\kappa_{NH}^{\Gamma}$. From Proposition~\ref{prop:inertial-mult}, we have
already seen that \(NH^{*,\Gamma}_T(T^*\C^n)\) can be written as
\[
NH^{*,\Gamma}_T(T^*\C^n) \cong \Q[u_1,u_2, \ldots,
u_n][\{\gamma_t\}_{t \in \Gamma}] \bigg/ {\mathcal I} +
{\mathcal J} + \left< \gamma_{id}-1 \right>.
\]
Thus it remains to describe each piece of the kernel, \({\mathcal K} =
\ker(\kappa_{NH}^{\Gamma}) = \bigoplus_{t \in \Gamma}
\ker(\kappa_{NH}^t).\)
Proposition~\ref{prop:kernel-kappa-t} implies that \(\gamma_t \star
{\mathcal K}_t \subseteq {\mathcal K}\) where ${\mathcal K}_t$ is
defined in~\eqref{eq:def-ideal-Kt} and is here
considered as an ideal in $NH_T^{*,id}(T^*\C^n)$. Note that \(\gamma_t \star {\mathcal
K}_t\) is  a subset of the $t$-th graded piece
\(NH^{*,t}_T(T^*\C^n)\).

This concludes the proof except for one subtlety: in
Proposition~\ref{prop:kernel-kappa-t} we described ideal generators
for ${\mathcal K}_t$ with respect to the {\em standard} ring structure
of Borel-equivariant cohomology, whereas in
Theorem~\ref{theorem:main-comb} we present generators with respect to
the $\star$ (or equivalently $\smile$) product on
$NH^{*,\Gamma}_T(T^*\C^n)$. Thus, given generators in the standard
product, it is not immediate that their union (multiplied by
appropriate $\gamma_t$) would yield ideal generators for ${\mathcal
K}$ in the $\star$ product. However, the $\id$-graded piece
$NH^{*,\id}_T(T^*\C^n) \cong H^*_T(pt)$ is a subring of
$NH^{*,\Gamma}_T(T^*\C^n)$ in the $\star$ product, and multiplication
in the $\star$ product of elements in $NH^{*,\id}_T(T^*\C^n)$ and
$NH^{*,t}_T(T^*\C^n)$ agrees with the standard $H^*_T(pt)$-module
structure on $H^*_T((T^*\C^n)^t)$ in Borel-equivariant
cohomology. Since Theorem~\ref{theorem:Hausel-Sturmfels} gives
$H^*_T(pt)$-module generators for each ${\mathcal K}_t$, the result
follows.
\end{proof}

\section{Examples}\label{sec:example}

We compute several explicit examples in this section to illustrate our
methods. Throughout, we identify Lie algebras with their dual spaces
using the standard inner product.
When illustrating the hyperplane arrangements, we will shade the
intersection of the positive half-spaces corresponding to the
cooriented hyperplanes.

\subsection{A hyperk\"ahler analogue of an orbifold $\P^2$}\label{subsec:orbifoldy-P2}

\begin{figure}[h]
\centerline{\epsfig{figure=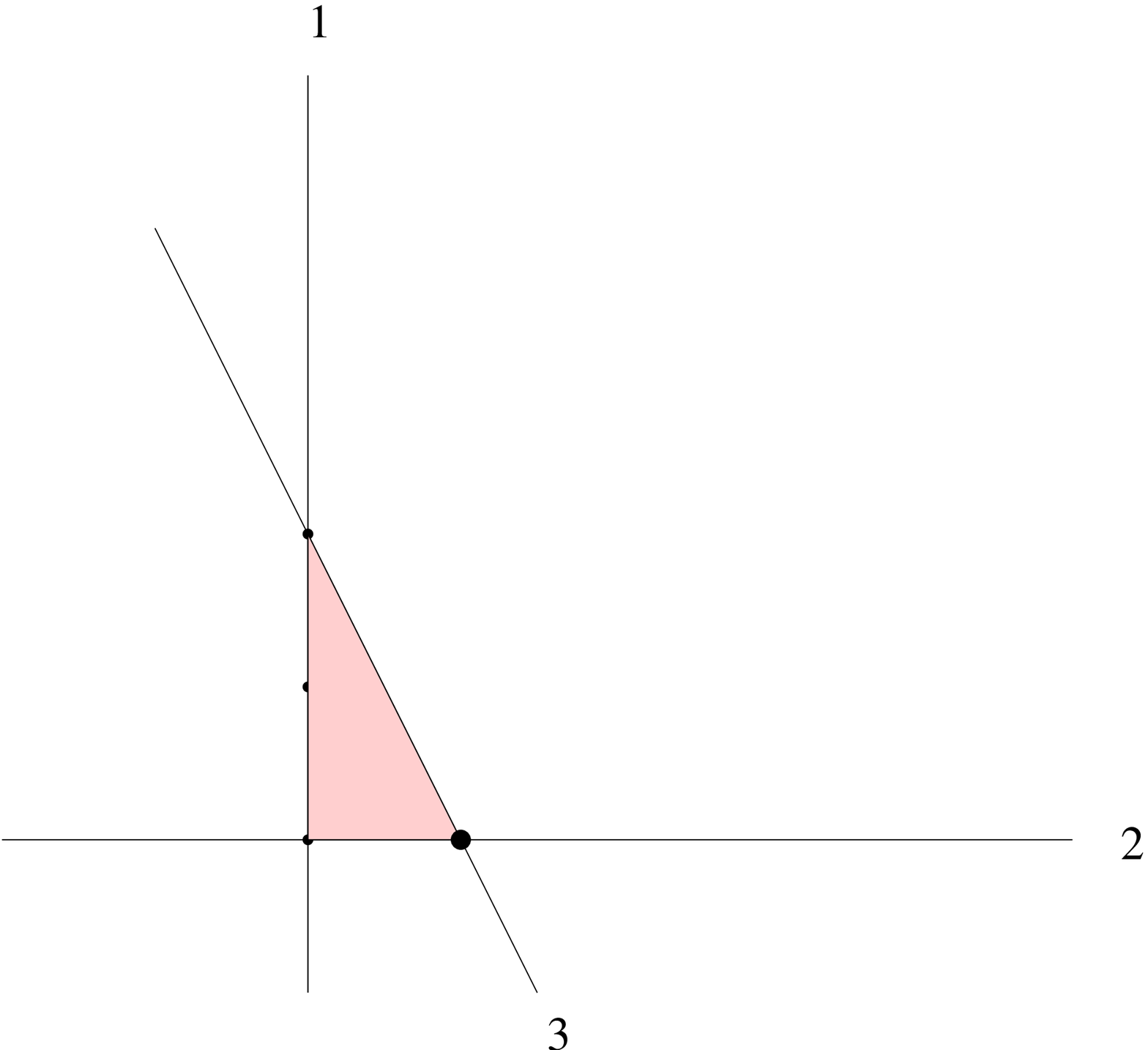, width=3in}}
\caption{An example of an orbifold hypertoric variety obtained by reducing
$\H^3$ by $S^1$. The corresponding K\"ahler toric variety is a $\P^2$
  with a single orbifold point, which maps to $H_2 \cap
  H_3$.}\label{fig:orbiP2}
\end{figure}

We begin with an example in which the corresponding K\"ahler toric
variety is an orbifold $\P^2$. Let ${\mathcal H}$ be the hyperplane
arrangement depicted in Figure~\ref{fig:orbiP2} and denote the
corresponding hypertoric variety by $M$. In this example,
\(n=3, d=2, k=1.\)
Here we will take the normal vectors to the hyperplanes to be primitive.

With respect to the standard bases in $\algt^3$ and $\algt^2$, the map $\beta$
in~\eqref{eq:Delzant-Lie} is given by
\begin{equation*}
\beta = \left[ \begin{array}{ccc} 1 & 0 & -2 \\ 0 & 1 & -1 \end{array} \right],
\end{equation*}
where the $i$-th column is the vector $a_i$ normal to the $i$-th
hyperplane in Figure~\ref{fig:orbiP2}. By Proposition~\ref{prop:finite-stabilizers}, the
single orbifold point maps to the intersection of the hyperplanes
$H_2$ and $H_3$.
The kernel of $\beta$ is given by the
span of the single vector $(2,1,1)$ in $\algt^3$.
Hence the $S^1$-action on $T^*\C^3$ with respect to which we
take a hyperk\"ahler quotient is induced by the linear action of $S^1$
on $\C^3$ with weights $2,1,1$ on the three coordinates,
respectively. In particular, it is immediate that the finite
stabilizer subgroup $\Gamma$ is just $\{\pm 1\} \cong \Z/2\Z$.
We compute the following table of logweights;
the quantity $2 \age(t)$
is the degree of the corresponding generator as in~\eqref{eq:def-age}.

\renewcommand{\arraystretch}{1.3}
\begin{equation}\label{eq:orbiP2-table}
\begin{array}{c|c|c|c|c|c|}
t & a_2(t) & a_1(t) & a_1(t) &  2 \ \age(t) & \genfrac{}{}{0pt}{0}{\mbox{generator of}}{NH_{S^1}^{*,t}(T^*\C^3)} \\
\hline \id & 0 & 0 & 0 &  0 & \gamma_{\id}  \\
\hline -1  & 0 & \frac{1}{2} & \frac{1}{2} & 4 & \gamma_{-1} \\ \hline
\end{array}
\end{equation}


Since there is only one non-trivial generator $\gamma_{-1}$ in
$NH_{S^1}^{*,\Gamma}(T^*\C^3)$ as a $H^*_{S^1}(pt)$-module, we only need to
compute a single relation of the form~\eqref{eq:inertial-relations},
namely, the product of $\gamma_{-1}$ with itself. Since \(t^2 = \id =
1\) for \(t = -1,\) we also have
\begin{equation*}
A(-1,-1)  = \{2,3\}, \qquad
B(-1,-1)  = \emptyset,\qquad
C(-1,-1)  = \emptyset,
\end{equation*}
as can be computed from the definitions~\eqref{eq:def-ABC},
and so we have
\[
\gamma_{-1}^2 -  u_2^2 u_3^2 \in {\mathcal I}.
\]
The ideal ${\mathcal J}$ of linear relations can be deduced from the matrix of
$\beta$ to be
\[
{\mathcal J} = \left< u_1 - 2u_3, u_2 - u_3 \right>.
\]
Finally, the ideal ${\mathcal K} = \ker(\kappa_{NH}^{\Gamma})$ may be
computed via the two pieces \(\ker(\kappa_{NH}^{\id}),
\ker(\kappa_{NH}^{-1}).\) We have from
Proposition~\ref{prop:kernel-kappa-t} that
\[
{\mathcal K}_{\id} = \left< u_1 u_2 u_3 \right>, \quad \mbox{and} \quad
{\mathcal K}_{-1} = \left< \gamma_{-1} u_1 \right>,
\]
from which we conclude that
\[
H^*_{CR}(M) \cong \Q[u_1,u_2,u_3,\gamma_{\id}, \gamma_{-1}] \bigg/
\left< \gamma_{-1}^2 - u_2^2 u_3^2, u_1 - 2u_3, u_2 - u_3, \atop u_1 u_2 u_3,
\gamma_{-1} u_1, \gamma_{\id} - 1 \right>,
\]
which is easily shown to be isomorphic to
\[
H^*_{CR}(M) \cong \Q[u, \gamma] \bigg/ \left< u^3, \gamma^2, \gamma u
\right>,
\]
where \(\deg(u) = 2, \deg(\gamma) = 4.\)
From this it is straightforward to compute that the orbifold
Poincar\'e polynomial for $M$ is given by
\[
P_{orb}(t, M) = 1 + t^2 + 2 t^4,
\]
so the orbifold Euler characteristic is $4$.

\begin{remark}
Let $M_n$ denote the hypertoric variety associated to the more general case in which $a_3 = (-n, -1)$ (so
the above case is $n=2$). The underlying K\"ahler toric variety is a
weighted $\P^2$ with a single orbifold point with orbifold structure group
$\Z/n\Z$. An analogous computation yields the orbifold Poincar\'e polynomial
\[
P_{orb}(t, M_n) = 1 + t^2 + n t^4,
\]
so $M_n$ has orbifold Euler characteristic $n+2$.
\end{remark}

\subsection{A quotient of $T^*\C^4$ by a $T^2$}\label{subsec:Nick}

\begin{figure}[h]
\centerline{\epsfig{figure=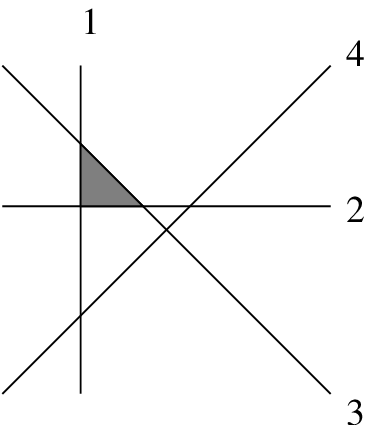}}
\caption{An example of an orbifold hypertoric variety obtained by reducing
$\H^4$ by $T^2$. The intersection $H_3 \cap H_4$ corresponds to the
  orbifold point.}\label{fig:orbiExample}
\end{figure}

We continue with an example in which the corresponding K\"ahler toric
variety is a smooth $\P^2$, but now we add an extra hyperplane which
introduces an orbifold point in the hypertoric variety.
Let ${\mathcal H}$ be the hyperplane arrangement depicted in
Figure~\ref{fig:orbiExample} and denote by $M$ the corresponding
hypertoric variety.
In this example, \(n=4, d=2, k=2.\)
We take primitive normals to these hyperplanes.


The map $\beta$ is given by the matrix
\begin{equation*}
\beta = \left[ \begin{array}{cccc} 1 & 0 & -1 & -1 \\ 0 & 1 & -1 & 1 \end{array} \right].
\end{equation*}
By
Proposition~\ref{prop:finite-stabilizers}, the single orbifold point
maps to the
intersection of the hyperplanes $H_3$ and $H_4$. The kernel of $\beta$ is given by the
Lie subalgebra \(\algt^k = \algt^2 \subseteq \algt^4\) given by the
span of the vectors \((1,1,1,0)\) and \((1,-1,0,1)\) in $\algt^4 \cong
\R^4$. Therefore, the $T$-action on $T^*\C^4$ with respect to which we
take a hyperk\"ahler quotient is given by
\begin{equation*}\label{eq:T2-action}
(t_1, t_2) \cdot (z, w) = 
(t_1 t_2 z_1, t_1 t_2^{-1} z_2, t_1 z_3, t_2 z_4, t_1^{-1} t_2^{-1} w_1, t_1^{-1} t_2 w_2, t_1^{-1}w_3, t_2^{-1}w_4),
\end{equation*}
where \(z =(z_1, z_2, z_3, z_4), w = (w_1, w_2, w_3, w_4),\) and here we have chosen
an identification of the kernel of $\exp(\beta)$ with the standard $2$-torus $T^2$.

The finite stabilizer group $\Gamma$ may now be computed as follows.
The weights
$\{\lambda_i\}_{i=1}^4$ for the $T^2$-action are given by \(\lambda_1
= (1,1), \lambda_2 = (1,-1), \lambda_3 = (1,0), \lambda_4 = (0,1).\)
The only minimal spanning subset which leads to a non-trivial
stabilizer is \(\{\lambda_1, \lambda_2\},\) and the stabilizer
subgroup is generated by the element \((-1, -1) \in T^2.\) Hence
\(\Gamma \cong \Z/2\Z.\) Then
immediately \(S(\id) = \{1,2,3,4\}\) and \(S\left((-1,-1)\right) = \{1,2\}.\) We will also
use the following table.

\renewcommand{\arraystretch}{1.3}
\begin{equation}\label{eq:orbiExample-table}
\begin{array}{c|c|c|c|c|c|c|}
t & a_{(1,1)}(t) & a_{(1,-1)}(t) & a_{(1,0)}(t) & a_{(0,1)}(t) &  2\ \age(t) & \genfrac{}{}{0pt}{0}{\mbox{generator of}}{NH_{T}^{*,t}(T^*\C^4)} \\
\hline \id & 0 & 0 & 0 & 0 & 0 & \gamma_{\id} \\
\hline (-1,-1) & 0 & 0 & \frac{1}{2} & \frac{1}{2} & 4 & \gamma_{(-1,-1)} \\ \hline
\end{array}
\end{equation}


As in the previous example, we only need to
compute a single relation of the form~\eqref{eq:inertial-relations},
namely, the product of $\gamma_{(-1,-1)}$ with itself.
We have
\begin{align*}
A((-1,-1),(-1,-1))  &= \{3,4\},\\
B((-1,-1),(-1,-1))  &= \emptyset,\\
C((-1,-1),(-1,-1))  &= \emptyset,
\end{align*}
and so
\[
\gamma_{(-1,-1)}^2  -  u_3^2 u_4^2 \in {\mathcal I}.
\]
The ideal of linear relations is
\[
{\mathcal J} = \left< u_1 - u_3 - u_4, u_2 - u_3 + u_4 \right>.
\]
Again as in the previous example, the ideal ${\mathcal K} =
\ker(\kappa_{NH}^{\Gamma})$ may be computed via the two pieces
\(\ker(\kappa_{NH}^{\id}), \ker(\kappa_{NH}^{(-1,-1)}).\) We have
\[
{\mathcal K}_{\id} = \left< u_2 u_3 u_4, u_1 u_3 u_4, u_1 u_2 u_4,
u_1 u_2 u_3 \right>, \quad \mbox{and} \quad
{\mathcal K}_{(-1,-1)} = \left< \gamma_{(-1,-1)} u_1, \gamma_{(-1,-1)} u_2 \right>.
\]
We conclude
\[
H^*_{CR}(M) \cong \Q[u_1,u_2,u_3,u_4,\gamma_{\id}, \gamma_{(-1,-1)}] /I
\]
where $$I=\left< \begin{array}{c} \gamma_{(-1,-1)}^2 - u_3^2 u_4^2, u_1 - u_3 - u_4, u_2 - u_3 +
u_4, \\
u_1 u_3 u_4, u_2 u_3 u_4, u_1u_2u_3, u_1u_2u_4\\
\gamma_{(-1,-1)} u_1, \gamma_{(-1,-1)}
u_2, \gamma_{\id} - 1
\end{array} \right>.$$
This simplifies to
\[
H^*_{CR}(M) \cong \Q[u_1, u_2, \gamma] \bigg/
\left< \gamma^2, u_1^3, u_2^3, u_1u_2^2, u_1^2 u_2, \gamma u_1, \gamma
u_2 \right>.
\]
Here, \(\deg(u_i) = 2, \deg(\gamma) = 4.\)
We see that the orbifold Poincar\'e
polynomial is
\[
P_{orb}(t, M) = 1 + 2t^2 + 4t^4,
\]
so the orbifold Euler characteristic is $7$.

\newcommand{\excise}[1]{}

\excise{
\subsection{An orbifold $\P^1 \times \P^2$}\label{subsec:P1-times-P2}

As a final example, we look at a hypertoric variety obtained by
taking a hyperk\"ahler quotient of $\H^5 \cong T^*\C^5$ by $T^2$. A
representation of the hyperplane arrangement (simplied for
visualization) is given in Figure~\ref{fig:P1timesP2}. Unlike the
previous examples, this orbifold hypertoric variety has nontrivial obstruction bundles.

\begin{figure}[h]
\psfrag{a1}{$a1$}
\psfrag{a2}{$a2$}
\psfrag{a3}{$a3$}
\psfrag{a4}{$a4$}
\psfrag{a5}{$a5$}
\centerline{\epsfig{figure=P1timesP2.eps, width=3in}}
\caption{An example of an affine hyperplane arrangement in
$\algt^3$ corresponding To an orbifold hypertoric variety obtained
by reducing $\H^5$ by $T^2$. It contains one $T^*\P^1$ subvariety with
a $\Z/2\Z$ global orbifold structure group, and another $T^*\P^1$
subvariety with a $\Z/3\Z$ global orbifold structure group. For the
sake of clear visuatlization, we have only depicted the compact
polytope corresponding to the corresponding orbifold K\"ahler toric
variety, but for the hypertoric case, the facets should be extended to
be hyperplanes.}\label{fig:P1timesP2}
\end{figure}

The Delzant exact sequence for this example is given by
\[
\xymatrix{
0 \ar[r] & \algt^2 \ar[r]^{\iota} & \algt^5 \ar[r]^{\beta} & \algt^3 \ar[r] & 0, \\
}
\]
where
\[
\beta = \left[ \begin{array}{ccccc} -3 & 0 & 1 & 0 & 0 \\
                                                             -2 & 1 & 0 & 0 & 0 \\
                                   0 & 0 & 0 & 1 & -1 \\
                                   \end{array} \right],
\]
with respect to the standard bases of $\algt^5$ and $\algt^3$. For a certain choice of basis for \(\ker(\beta) \cong \algt^2,\) we get that the dual map $\iota^*$ is given by
\[
\iota^* = \left[ \begin{array}{ccccc} 1 & 2 & 3 & 0 & 0 \\ 0 & 0 & 0 & 1 & 1 \end{array}\right],
\]
so the $T^2$-weights are given by \(\lambda_1 = (1,0), \lambda_2 =
(2,0), \lambda_3 = (3,0), \lambda_4 = (0,1), \lambda_5 = (0,1) \in
(\algt^2)^*_{\Z}.\) It is straightforward to see that \(\Gamma =
\left< \zeta \right>\) where \(\zeta = (e^{2\pi i/6}, 1) \in T^2 \cong
(S^1)^2.\) So \(\Gamma \cong \Z/6\Z.\)
We compute
\[
S(\id) = S((1,1)) = \{1,2,3,4,5\}, \quad S(\zeta) = S(\zeta^5) = \{4,5\}, \quad
S(\zeta^2) = S(\zeta^4) = \{3,4,5\}, \quad S(\zeta^3) = \{2,4,5\},
\]
and we also have the following table of logweights.

\renewcommand{\arraystretch}{1.3}
\begin{equation}\label{eq:orbiP2table}
\begin{array}{c|c|c|c|c|c|c|c|}
t & a_{(1,0)}(t) & a_{(2,0)}(t) & a_{(3,0)}(t) & a_{(0,1)}(t) & a_{(0,1)}(t) & 2 \ \age(t) & \genfrac{}{}{0pt}{0}{\mbox{generator of}}{NH_{T^2}^{*,t}(T^*\C^5)} \\
\hline \id & 0 & 0 & 0 & 0 & 0 & 0 & \gamma_{\id} \\
\hline \zeta & \frac{1}{6} & \frac{1}{3} & \frac{1}{2} & 0 & 0 & 6 & \gamma_1 \\
\hline \zeta^2 & \frac{1}{3} & \frac{2}{3} & 0 & 0 & 0 & 4 & \gamma_2 \\
\hline \zeta^3 & \frac{1}{2} & 0 & \frac{1}{2} & 0 & 0 & 4 & \gamma_3 \\
\hline \zeta^4 & \frac{2}{3} & \frac{1}{3} & 0 & 0 & 0 & 4 & \gamma_4 \\
\hline \zeta^5 & \frac{5}{6} & \frac{2}{3} & \frac{1}{2} & 0 & 0 & 6 & \gamma_5 \\ \hline
\end{array}
\end{equation}

\medskip
As an illustration, we now explicitly compute $\gamma_5 \star
\gamma_5$.  Here, \(t_1 = t_2 = \zeta^5,\) and we have \(S(t_1)^c =
S(t_2)^c = \{1,2,3\} = S(t_1)^c \cap S(t_2)^c.\) Using the
definitions of the subsets $A(t_1,t_2), B(t_1,t_2), C(t_1,t_2)$ as
given in Section~\ref{sec:comb}, we find that
\begin{equation}
A(t_1, t_2)  = \{3,4,5\}, \qquad
B(t_1, t_2)  = \emptyset, \qquad
C(t_1, t_2)  = \{1,2\}.
\end{equation}
Hence we may conclude, using Proposition~\ref{prop:inertial-mult}, that
\[
\gamma_5 \star \gamma_5 = - u_1 u_2 u_3^2 u_4^2 u_5^2 \cdot \gamma_4.
\]
The rest of the multiplication table may be computed similarly; we summarize the results below.
\medskip

\begin{equation}\label{eq:mult-table-orbiP2}
\begin{array}{c||c|c|c|c|c|}
\star & \gamma_1 & \gamma_2 & \gamma_3 & \gamma_4 & \gamma_5 \\ \hline \hline
\gamma_1  &  -u_1 u_2 u_3^2 \gamma_2 & u_1 u_2^2 \gamma_3 & u_1 u_3^2 \gamma_4 & u_1 u_2 \gamma_5 & -u_1^2 u_2^2 u_3^2 \gamma_{\id}  \\ \hline
\gamma_2  & & -u_1 u_2 \gamma_4  & -u_1 \gamma_5 & u_1^2 u_2^2
\gamma_{\id} & u_1 u_2 \gamma_1 \\ \hline
\gamma_3 & & & u_1^2 u_3^2 \gamma_{\id} & u_1 \gamma_1 & -u_1 u_3^2 \gamma_2 \\ \hline
\gamma_4 & & & & -u_1 u_2 \gamma_2 & -u_1 u_2^2 \gamma_3 \\ \hline
\gamma_5 & & & & & -u_1 u_2 u_3^2 \gamma_4 \\ \hline
\end{array}.
\end{equation}

Thus, as a ring, the inertial cohomology is
\[
\NHGTCot \cong \Q[u_1, u_2, u_3, u_4, u_5, \gamma_{\id} \gamma_1, \gamma_2, \gamma_3, \gamma_4, \gamma_5] \big/ {\mathcal I} +
{\mathcal J} + \left< \gamma_{\id} - 1 \right>,
\]
where ${\mathcal I}$ is the ideal generated by the multiplication relations~\eqref{eq:mult-table-orbiP2},
and ${\mathcal J}$ is the ideal $\left< \im(\beta^*)\right>$, which in this case is
\begin{equation}\label{eq:mathcalJ-ex2}
{\mathcal J} = \left< -3u_1+u_3, -2u_1+u_2, u_4- u_5 \right>.
\end{equation}
We now need to compute the kernel of $\kappa^t_{NH}$ for each $t \in
\Gamma$. Since \((T^*\C^5)^{\zeta^2} =
(T^*\C^5)^{\zeta^4},\) however, we only need to compute three separate
(ordinary) Kirwan kernels.

For the case \(t = \id,\) this is just the usual computation of the
ordinary cohomology of a hypertoric variety, and Theorem~\ref{theorem:Hausel-Sturmfels} immediately yields
\[
{\mathcal K}_{\id} = \left< u_1 u_2 u_3, u_4 u_5 \right>.
\]
For the others, we demonstrate the computation for the case \(t =
\zeta^2.\). Since \(S(\zeta^2)^c = \{1,2\},\)
the ambient subspace of the affine hyperplane arrangement \({\mathcal
H}_{\zeta^2}\) is the $1$-dimensional intersection \(H_1 \cap H_2,\)
as indicated in Figure~\ref{fig:gamma2-subarr}. This is the moment
image of the hypertoric subvariety $M_{\{3,4,5\}} \subseteq M$ with global
orbifold structure group
\[
\mathrm{span}_{\Q}\{a_1, a_2\} \cap \algt^3_{\Z} \bigg/ \mathrm{span}_{\Z}\{a_1, a_2\} \cong \Z/3\Z.
\]
Since \(H_3 \cap (H_1 \cap H_2) = \emptyset,\) this corresponds to an ``extra'' generator as in Section~\ref{sec:comb}, while $H_4$ and $H_5$ both intersect nontrivially with \(H_1 \cap H_2\) and yield $0$-dimensional affine hyperplanes (i.e. points) in $H_1 \cap H_2$. See Figure~\ref{fig:gamma2-subarr}.

\begin{figure}[h]
\centerline{\epsfig{figure=gamma2-subarr.eps, width=3in}}
\caption{The hyperplane arrangement ${\mathcal H}_{\zeta^2}$ is
  obtained by intersecting the hyperplanes $H_3, H_4,$ and $H_5$ with
  $H_1 \cap H_2$, where the latter is indicated here by the bold vertical line. This
  intersection $H_1 \cap H_2$ is exactly the moment image of the
  subvariety $M_{\{3,4,5\}}$ with global orbifold structure group
  $\Z/3\Z$. The intersections of the hyperplanes $H_4$ and $H_5$
  produce the hyperplanes ($0$-dimensional, in this case, so just
  points) in $H_1 \cap H_2$ indicated by the bold dots. The hyperplane
  $H_3$, indicated here by a shaded region, does not intersect $H_1
  \cap H_2$, hence corresponds to an ``extra'' generator as in
  Proposition~\ref{prop:kernel-kappa-t}.}\label{fig:gamma2-subarr}
\end{figure}

Applying Proposition~\ref{prop:kernel-kappa-t} immediately yields the relations $\left< u_3, u_4 u_5 \right>$, and hence we conclude
\[
{\mathcal K}_{\zeta^2} = \left< \gamma_2 u_3, \gamma_2 u_4 u_5 \right>.
\]
Proceeding similarly for the other \(t \in \Gamma,\) we finally obtain
\[
H^*_{CR}(M) \cong
\Q[u_1, u_2, u_3, u_4, u_5, \gamma_1, \gamma_2, \gamma, \gamma_4, \gamma_5] \bigg/
{\mathcal I} + {\mathcal J} +
\left< u_4 u_5, u_1 u_2 u_3, \gamma_1, \gamma_2 u_3,  \gamma_3 u_2 , \gamma_4
u_3, \gamma_5 \right>,
\]
where ${\mathcal I}$ is generated by the
relations~\eqref{eq:mult-table-orbiP2} and ${\mathcal J}$ is given in~\eqref{eq:mathcalJ-ex2}.
}

\section{Appendix: Inertial cohomology and Chen-Ruan cohomology}\label{sec:definitionChase}

In this section, we show that there is a natural equivalence between
the inertial cohomology of a stably complex space $Z$ from
Section~\ref{sec:definition} and the orbifold cohomology of $Z/T$ when
$T$ acts locally freely, i.e. that there exists a graded ring
isomorphism
\begin{equation}\label{eq:NHisomCR}
NH_T^{*,\diamond}(Z) \cong H_{CR}^*(Z/T).
\end{equation}
Applying this isomorphism to the case when $Z$ is a level
set of the hyperk\"ahler moment map on $T^*\C^n$ and $M=Z/T$ is an
orbifold hypertoric variety
completes the proof of
Theorem~\ref{th:surjectiontoCRcohomology}. (There is a proof of a
similar statement in \cite{GHK05}, but here we drop their compactness assumption.) In addition, we show that the
definition of the product structure for $NH^{*,\diamond}_T(Z)$
is also equivalent to another description used in the algebraic-geometry literature
(e.g. \cite{FG03}, \cite{BCS05}).

We first prove~\eqref{eq:NHisomCR} as additive groups. We simplify
the presentation in \cite{CR04} to the case when the group involved is
abelian, and $X=Z/T$. Let $X_t := \{(p,t): p\in Z/T, t\in G_p\}$,
where $G_p$ is the local orbifold structure group at the point $p\in Z/T$. We assume
for simplicity that $X_t$ is connected (if not, take a direct sum over connected
components).  By definition,
\begin{equation}\label{eq:def-CR-cohomology}
H_{CR}^d(X) := \bigoplus_{t\in T} H^{d-2\sigma_t}(X_t),
\end{equation}
where the degree shift $\sigma_t$ is constant on connected components,
and is defined below.  Since in our case \(X = Z/T\) is a global
quotient, each $G_p$ is a subgroup of $T$ and $X_t= \{(zT,t): zT\in
Z/T,\ t\in Stab(z)\}.$ In other words, $X_t = Z_t/T,$ where
$Z_t := \{(z,t): z\in Z,\ t\in Stab(z)\}$ and the $T$-action
is on the first coordinate. Notice that, for a given $t\in T$, $Z_t
\iso Z^t := \{z\in Z: t\cdot z = z\}$. Since $T$ acts locally freely
on $Z$, it certainly acts locally freely on $Z^t$. Therefore,
$H^*(Z^t/T) \cong H^*_T(Z^t)$ (with $\Q$ coefficients) and
\begin{equation}\label{eq:CR-cohom-as-eqvt-cohomology}
H_{CR}^d(Z/T) \cong \bigoplus_{t\in T} H_T^{d-2\sigma_t}(Z^t)
\end{equation}
as additive groups. The right hand side
of~\eqref{eq:CR-cohom-as-eqvt-cohomology} is exactly the definition of
$NH^{*,\diamond}_T(Z)$, so we have proved the
additive isomorphism~\eqref{eq:NHisomCR}.

We now prove that the isomorphism~\eqref{eq:NHisomCR} also preserves
the grading. The number $\sigma_t$ appearing
in~\eqref{eq:def-CR-cohomology} is obtained as
follows; we assume \(X_t \neq \emptyset.\)
At any point \(p \in X,\) let \(\rho_p: G_p \to GL(k,\C)\) be a
representation specifying a local model \(\C^k/G_p\) at $p$. Since
$G_p$ is abelian, the image of $\rho_p$ is simultaneously
diagonalizable; denote by \(\{a_{\lambda_j}(t)\}_{j=1}^k\) the
logweights of the eigenvalues of \(\rho_p(t), t \in G_p.\)
The sum
\begin{equation}\label{eq:def-sigmat}
\sigma_t := \sum_{j=1}^k a_{\lambda_j}(t) \in \Q
\end{equation}
is well-defined, constant on connected components of $X_t$, and gives
the degree shift in \cite{CR04} and~\eqref{eq:def-CR-cohomology}.


We now show that this degree shift encoded by $\sigma_t$ agrees with
the degree shift in the definition of the grading for inertial
cohomology in \cite{GHK05}. The local model \(\C^k/G_p\)
can also be obtained by looking at the original $T$-space $Z$. Namely,
given a lift $z$ of the point \(p \in X,\) $G_p$ is exactly $\Stab(z)
\subseteq T$ and the representation $\rho_p$ above is given by the
action of $G_p$ on the normal bundle \(\nu(T \cdot z, Z)\) in
$T_zZ$. Moreover, since $t$ acts trivially on $Z^t$, the only
nontrivial eigenvalues of $\rho_p(t)$ are those which occur in the
representation of $\langle t \rangle$ on a further quotient
$\nu(Z^t, Z)$. In particular one may conclude that the
sum~\eqref{eq:def-sigmat} equals the sum
\begin{equation}\label{eq:def-age}
age(t) := \sum_{\C_{\lambda}\subset \nu(Z^t,Z)} a_{\lambda}(t).
\end{equation}
Even if $X$ is not compact, the grading shift is well defined (as long as $X$ is finite dimensional). In particular, the normal bundle $\nu(Z^t,Z)$ does not degenerate as it goes out to infinity. This shows that the gradings
agree.

We have left to show that the isomorphism~\eqref{eq:NHisomCR}
preserves the ring structure. The products on both
$NH^{*,\diamond}_T(Z)$ and $H^*_{CR}(X)$ are defined using the notion of
an obstruction bundle, so we begin by showing that the obstruction
bundle of Definition~\ref{def:obstructionbundle}, defined upstairs on $Z$,
descends to the obstruction bundle of Chen and Ruan, defined on the
quotient \(X = Z/T.\) In their original paper~\cite{CR04}, the authors
define these bundles over $3$-twisted sectors; however, their
construction can be greatly simplified in the case of a global
quotient \(X = Z/T,\) so we restrict attention to this case below.

Chen and Ruan define their obstruction bundle using two ingredients;
we describe each in turn. Consider a point \([z] \in X = Z/T,\) and
suppose that \(t_1, t_2 \in G_{[z]} \subseteq T.\) Let \(H := \langle
t_1, t_2 \rangle\) be the finite subgroup they generate. Let
\[
X_{(t_1, t_2, t_3)} := \{(p,t_1,t_2,t_3): p \in X, \hs \hs t_1, t_2,
t_3 \in G_p, \hs t_1 t_2 t_3 = 1 \}.
\]
Then there is a smooth map \(e: X_{(t_1,t_2,t_3)} \to X\) projecting
to the first term. Let \(e^*TX\) be the pullback of the tangent
bundle; this is a complex $H$-equivariant orbi-bundle over
$X_{(t_1,t_2, t_3)}$ and is the first ingredient in the Chen-Ruan
definition of the obstruction bundle.


The second ingredient involves only the subgroup $H$. Let $\Sigma= (\Sigma,t_1,t_2,t_3,H)$ be
a proper smooth Galois $H$ cover of $\P^1$ branched over
$\{0,1,\infty\}$ (for details see \cite{FG03}[Appendix]). The
$H$-action on $\Sigma$ induces an $H$-action on $H^1(\Sigma,{\mathcal
  O}_{\Sigma})$, so we may define the topologically trivial
$H$-equivariant bundle with fiber $H^1(\Sigma,{\mathcal O}_{\Sigma})$
over $X_{(t_1,t_2,t_3)}$ of complex rank $\mbox{genus}(\Sigma)$,
where the $H$-action is only on the fiber. We denote this bundle by
$H^1(\Sigma,{\mathcal O}_{\Sigma})$.  Then the obstruction bundle of
Chen and Ruan is given by the $H$-invariant part of the tensor product
of these two bundles, i.e.
\begin{equation}\label{eq:def-CR-obstruction}
E := (H^1(\Sigma, {\mathcal O}_{\Sigma}) \otimes e^*TX)^H.
\end{equation}


We now wish to show that the obstruction bundle of
Definition~\ref{def:obstructionbundle} descends
to~\eqref{eq:def-CR-obstruction}. As a first step, observe that
$X_{(t_1,t_2,t_3)}$ is isomorphic to $Z^H/T$, so the base spaces of
the two bundles certainly correspond. One reasonable way to lift the bundle might be to replace $e^*TX$ with $e^*TZ$ in the Chen-Ruan definition. However, this tangent bundle is not complex. Since a fiber $(e^*TX)_{[z]}$ of the orbi-bundle $e^*TX$ can be constructed via
$T$-equivalence classes in $\nu(T \cdot z, Z)$, a natural idea would be to split $e^*TZ$ at any point $z$ into the tangent directions along the orbits (which should not contribute), and its (complex) quotient bundle, $\nu(T\cdot z,Z)$.  Alternatively, one can split $e^*TZ$ into the tangent directions $TZ^H$ along the fixed point set, and its (complex) quotient $\nu(Z^H,Z)$. In either case, $TZ^H$ (or its quotient in $\nu(T\cdot z, Z)$) does not contribute to the obstruction bundle, since
\begin{equation}\label{eq:tangent-bundle-trivial}
(H^1(\Sigma,\mathcal{O}_{\Sigma})\otimes TZ^H)^H =
 H^1(\Sigma,\mathcal{O}_{\Sigma})^H\otimes TZ^H = H^1(\C
 P^1,\mathcal{O}_{\C P^1})\otimes TZ^H=0.
\end{equation}
Thus only the normal bundle $\nu(Z^H,Z)$ contributes, and we see that
\begin{equation}\label{eq:lift-of-CR-bundle}
\widetilde{E} := (H^1(\Sigma, {\mathcal O}_{\Sigma}) \otimes \nu(Z^H,
Z))^H \to Z^H
\end{equation}
quotients to $E$. Note that $\nu*(Z^H,Z)$ is well-defined, even if $Z^H$ is not compact. $Z^H$ is a closed submanifold containing as a submanifold the orbit through $z$; thus $\nu(Z^H,Z)$ is a quotient of a local model on $\nu(T\cdot z,Z)$ of the representation downstairs.

It remains to show that the only $H$-invariant subspaces of
$\nu(Z^H,Z)$ which contribute to~\eqref{eq:lift-of-CR-bundle} are the
$H$-isotypic components $\C_{\lambda} \subseteq I_{\lambda} \subseteq \nu(Z^H, Z)$ with
\(a_{\lambda}(t_1)+a_{\lambda}(t_2)+a_{\lambda}(t_3) = 2.\) We analyze
each piece $\C_{\lambda}$ separately. We use \v{C}ech cohomology to
compute with the $H^1(\Sigma,{\mathcal O}_{\Sigma})$, so let
$\mathcal{U} = \{U_i\}_{i \in I}$ be an $H$-invariant open cover of
$\Sigma$, i.e. for every \(\sigma \subseteq I\) there exists $\tau$
such that \(h \cdot U_{\sigma} = U_{\tau}.\) We denote this
by \(h \cdot \sigma = \tau\) for simplicity.

We claim that for \(z \in Z^H,\) the fiber $(\C_{\lambda} \otimes
H^1(\Sigma,{\mathcal O}_{\Sigma}))^H$ is isomorphic as a $\C$-vector
space to $H^1(\Sigma,{\mathcal L}_{\lambda})$, where ${\mathcal
  L}_{\lambda}$ is the sheaf of $H$-invariant sections of the
topologically trivial $H$-equivariant line bundle $L_{\lambda} =
\Sigma \times \C_{\lambda}$ over $\Sigma$. This can be seen at the
level of cochains by the map
\[
\xymatrix @R=-0.2pc {
\phi: (\C_{\lambda} \otimes C^{\ell}({\mathcal U}, {\mathcal
  O}_{\Sigma}))^H \ar[r] & C^{\ell}({\mathcal U}, {\mathcal
  L}_{\lambda}) \\
z \otimes s \ar@{|->}[r] & zs. \\
}
\]
It is straightforward to check that $\phi$ is well-defined and an
isomorphism using the definition of the $H$-action on $\C_{\lambda}
\otimes C^{\ell}({\mathcal U}, {\mathcal O}_{\Sigma})$, which can be
written, for $\sigma \subseteq I$,
\[
\left( h \cdot (z \otimes s) \right) \vert_{\sigma} = e^{2\pi
  ia_{\lambda}(h)} z \otimes s \vert_{h^{-1} \sigma}.
\]
It may also be checked that $\phi$ commutes with the \v{C}ech
differential, so \((\C_{\lambda} \otimes H^1(\Sigma, {\mathcal
  O}_{\Sigma}))^H \cong H^1(\Sigma, {\mathcal L}_{\lambda}),\) as
desired.

Furthermore, it is shown in \cite{BCS05} that \(H^1(\Sigma, {\mathcal
  L}_{\lambda}) \cong H^1(\P^1, {\mathcal O}(-a_{\lambda}(t_1) -
  a_{\lambda}(t_2) - a_{\lambda}(t_3))).\) The latter is
  $1$-dimensional exactly when the sum inside is $-2$ and $0$-dimensional
  otherwise, so $(\C_{\lambda} \otimes H^1(\Sigma, {\mathcal
  O}_{\Sigma}))^H$ contributes nontrivially to $\widetilde{E}$ if and
  only if \(a_{\lambda}(t_1) + a_{\lambda}(t_2) + a_{\lambda}(t_3) =
  2.\) As a bundle, each of these contributions is a line bundle over
  $Z^H$ given as a sub-bundle of $\nu(Z^H, Z)$, since by construction
  $H^1(\Sigma, {\mathcal O}_{\Sigma})$ is the trivial bundle over
  $Z^H$. We conclude that
\[
\widetilde{E} \cong \sum_{I_{\lambda} \subseteq \nu(Z^H,Z) \atop
  a_{\lambda}(t_1)+a_{\lambda}(t_2) + a_{\lambda}(t_3) = 2} I_{\lambda},
\]
where $I_\lambda$ is the isotypic component of $\nu(Z^H,Z)$ of weight $\lambda$.

Finally, under the isomorphism $H_T^*(Z^H)\iso
H^*(Z^H/T)$, the equivariant Euler class $e_T(\widetilde{E})$ is mapped to the
ordinary Euler class $e(E)\in H^*(Z^H/T)$. The
$\smile$-product is then constructed to be identical to the definition
given in \cite{CR04}.  We have proven~\eqref{eq:NHisomCR}, which we
record as follows.
\begin{theorem}
The inertial cohomology $NH^{*,\diamond}_{T}(Z)$ is isomorphic as a graded ring to the orbifold cohomology $H^*_{CR}(M)$.
\end{theorem}

We now prove the correspondence of our
definition of
the obstruction bundle with a description in terms of right derived
functors used
in the algebraic geometry literature (e.g. \cite{BCS05},
\cite{FG03}). For this exercise, it is convenient to use the
description in~\eqref{eq:lift-of-CR-bundle}. In the
algebraic-geometric context, the definition of the obstruction bundle
(in the case of a global quotient by a locally free action) over $Z^H$
is given as \(R^1\pi^H_*(\pi^*TZ|_{Z^H}),\) where $\Sigma, H$ are as
above, \(\pi: Z^H \times \Sigma \to Z^H\) is the projection,
$\pi_*^H$ is the functor
``pushforward and take $H$-invariants'', and $R^1\pi_*^H$ is its first
right derived functor. By an argument similar
to~\eqref{eq:tangent-bundle-trivial}, only the normal bundle
$\nu(Z^H,Z)$ contributes nontrivially to this bundle, so
\(R^1\pi_*^H(\pi^*TZ|_{Z^H}) = R^1\pi_*^H(\pi^*\nu(Z^H,Z)).\) We will
work with this second description; in particular, we will show that
in our (not necessarily algebraic) context, the right hand side
of this equation is equal to our
bundle~\eqref{eq:lift-of-CR-bundle}.

We begin by computing $R^1\pi_*(\pi^*\nu(Z^H,Z))$. The sheaf of
sections of $\pi^*\nu(Z^H,Z)$ is a ${\mathcal O}_{Z^H \times
  \Sigma}$-module, where ${\mathcal O}_{Z^H \times \Sigma}$ is the
sheaf of smooth functions on $Z^H \times \Sigma$ that are holomorphic
restricted to any fiber of $\pi$. By the push-pull formula,
\[
R^1\pi_*(\pi^*\nu(Z^H,Z)) = \nu(Z^H,Z) \otimes R^1\pi_*({\mathcal
  O}_{Z^H \times \Sigma}).
\]
Moreover, the pushforward sheaf $\pi_*({\mathcal O}_{Z^H \times
  \Sigma})$ can be described as
\[
\pi_*({\mathcal O}_{Z^H \times \Sigma}) = {\mathcal O}_{Z^H} \otimes
\Gamma(\Sigma,{\mathcal O}_{\Sigma}),
\]
where here ${\mathcal O}_{Z^H}$ is the sheaf of smooth functions on $Z^H$
and ${\mathcal O}_{\Sigma}$ is the (usual) sheaf of holomorphic
functions on $\Sigma$. This implies that $R^1\pi_*({\mathcal O}_{Z^H
  \times \Sigma}) = {\mathcal O}_{Z^H} \otimes H^1(\Sigma, {\mathcal
  O}_{\Sigma})$, so we finally have
\[
R^1\pi^H_*(\pi^*\nu(Z^H,Z)) = \left(\nu(Z^H,Z) \otimes H^1(\Sigma,
{\mathcal O}_{\Sigma})\right)^H \to Z^H,
\]
as desired.

\section*{Acknowledgments}
It is our pleasure to thank the
American Institute of Mathematics for hosting a conference on the
subject of Kirwan surjectivity, at which the authors began work on
this project. We also thank Nicholas Proudfoot
for useful conversations. The first author thanks Lisa Jeffrey and the University
of Toronto for hospitality while some of this work was being
conducted. The second author similarly thanks George Mason University.
RG was partially supported by NSF-DMS Grant 0305128.

\end{document}